%% file: amsucnotes01.tex
\documentclass{conm-p-l}
\usepackage{amsmath,amsthm}
\usepackage{amssymb}
\usepackage[all]{xypic}
\usepackage{amsbsy}

\newtheorem{theorem}{Theorem}[section]
\newtheorem{lemma}[theorem]{Lemma}

\newtheorem{prop}[theorem]{Proposition}
\newtheorem{proposition}[theorem]{Proposition}
\newtheorem{corollary}[theorem]{Corollary}

\theoremstyle{definition}
\newtheorem{definition}[theorem]{Definition}
\newtheorem{example}[theorem]{Example}

\newtheorem{examples}[theorem]{Examples}
\newtheorem{proof-of-halfcm5}[theorem]{Proof of Proposition \ref{half-cm5}}

\theoremstyle{remark}
\newtheorem{remark}[theorem]{Remark}
\newtheorem{remarks}[theorem]{Remark}

\numberwithin{equation}{section}
\input{macros}

\begin{document}

\title{Model Categories and Simplicial Methods}

\author{Paul Goerss}
\address{Department of Mathematics, Northwestern University, Evanston IL 60208}
\email{pgoerss@math.northwestern.edu}
\thanks{The first author was partially supported by the National Science Foundation.}

\author{Kristen Schemmerhorn}
\address{Department of Mathematics \& Computer Science, Albion College, Albion, MI 49224}
\email{kschemmerhorn@dom.edu}
\curraddr{Dominican University, River Forest, IL 60305}

\thanks{Many thanks to Srikanth Iyengar, Ieke Moerdijk, and Brooke
Shipley, all of whom went through this manuscript very carefully and made
many helpful suggestions.}

\date{}

\begin{abstract}
There are many ways to present model categories, each with a 
different point of view.  Here we'd like to treat model categories as a
way to build and control resolutions. This an historical approach, as in
his original and spectacular applications of model categories, Quillen
used this technology as a way to construct resolutions
in non-abelian settings; for example, in his work on the homology of
commutative algebras \cite{AQ}, it was important to be very flexible
with the notion of a free resolution of a commutative algebra. Similar
issues arose in the paper on rational homotopy theory \cite{rational}.
(This paper is the first place where the now-traditional
axioms of a model category are enunciated.) We're going to
emphasize the analog of projective resolutions, simply because these
are the sort of resolutions most people see first. Of course, the theory is
completely flexible and can work with injective resolutions as well.

There are now any number of excellent sources for getting into the subject
and since this monograph is not intended to be complete, perhaps the
reader should have some of these nearby. For example,
the paper of Dwyer and Spalinski \cite{DS} is a superb and short introduction,
and the books of Hovey \cite{Model} and Hirschhorn \cite{Local} provide
much more in-depth analysis. For a focus on simplicial model categories
-- model categories enriched over simplicial sets in an appropriate way
-- one can read \cite{GJ}. Reaching back a bit further, there's no
harm in reading the classics, and Quillen's
original monograph \cite{HA} certainly falls into that category.
\end{abstract}

\maketitle

\tableofcontents

\section{Model Categories and Resolutions}

\subsection{Chain complexes}

Let us begin with a familiar and basic example. Nothing in this section
is supposed to be new, except possibly the point of view.

Let $R$ be a commutative ring. Then we will denote the category of
$R$-modules by $\Mod_{R}$. In this category, there is a
distinguished class of objects in $\Mod_{R}$; namely,
the projective $R$-modules. There is also  a distinguished class of morphisms
in $\Mod_{R}$, the surjections. These two
classes determine each other; indeed, an $R$-module $P$ is projective
if and only if
$$
\Mod_{R}(P,f): \Mod_R(P,M) \to \Mod_R(P,N)
$$
is a surjective map of sets for all surjections $f:M \to N$ of $R$-modules. 
Conversely, a morphism $f$ is a surjection if and only if
$\Mod_{R}(P,f)$ is surjective for all projectives $P$.  

There is a distinguished projective in $\Mod_R$, namely $R$ itself, and
we can use that fact to show that $\Mod_{R}$ has enough
projectives. This means that for all $M \in \Mod_R$, there is a
surjection $P \to M$ with $P$ projective. (Here we combine the notions
of a projective module and a surjective morphism.)
If $M \in \Mod_{R}$, we define
$$
P(M) = \oplus_{f:R \to M} R
$$
where $f$ runs over the morphims in $\Mod_R$ from $R$ to $M$; then evaluation
defines a surjection $\epsilon_M:P(M) \to M$. We have written
this morphim as if it were a functor of $M$, and indeed it is. This will
be useful later.

Now let $\mathbf{Ch}_{\ast}(R)$ be the chain complexes in $\Mod_{R}$.
By this we mean the non-negatively graded chain complexes, at least for
now. These have the form
\[
\cdots \longrightarrow C_{3} \longrightarrow C_{2} \longrightarrow
C_{1} \longrightarrow C_{0}.
\]
Then a projective resolution of a module $M \in \Mod_{R}$, is
a chain complex $P_{\bullet}$ in $\mathbf{Ch}_{\ast}R$, so that
\begin{enumerate}

\item each $P_n$ is a projective;

\item $H_n P_\bullet = 0$ for $n > 0$; and

\item there is a morphism $P_0$ to $M$ which induces an isomorphim
$H_0 P_\bullet \cong M$.
\end{enumerate}

We can rephrase the last two points. We regard $M$ as a chain complex
concentrated in degree $0$ and we can say that there is a morphism of chain
complexes $P_\bullet \to M$ which induces an isomorphism on homology.

Chain complexes of projectives have the following important property. 
Recall that a chain complex $N_\bullet$ is {\it acyclic} if 
$H_n N_\bullet = 0$ for $n > 0$. Then there is automatically a
homology isomorphism $N_\bullet \to H_0N_\bullet$, where the target is
regarded as the chain complex concentrated in degree zero. Now suppose
$P_\bullet$ is a chain complex so that each $P_n$ is projective; we are
not assuming it is acyclic. Then, given any morphism of chain complexes
$f:P_\bullet \to H_0N_\bullet$ (or, equivalently, a morphism
$H_0P_\bullet \to H_0N_\bullet$ of modules)
we can solve the lifting problem
\begin{equation}\label{proto-lifting}
\xymatrix{
& N_{\bullet} \ar[d]^-{\epsilon}\\
P_{\bullet} \ar@{-->}[ru]^-{g} \rto_-f& H_0N_\bullet\\
}
\end{equation}
and $g$ is unique up to chain homotopy. This fact is then used to
prove the uniqueness of projective resolutions up to chain
homotopy.

It turns out that we can solve a much more general lifting problem
than in (\ref{proto-lifting}). Here is the result.

\begin{prop}\label{half-cm5} Suppose we have two morphisms
$j:A_\bullet \to B_\bullet$ and $q: M_{\bullet} \to N_{\bullet}$
of chain complexes so that
\begin{enumerate}

\item for all $n \geq 0$, $A_n \to B_n$ is an injection and
$B_n/A_n$ is projective; and

\item $H_\ast q$ is an isomorphism and $M_{n} \longrightarrow N_{n}$ is a
surjection for $n > 0$.
\end{enumerate}
Then any lifting problem
\[
\xymatrix{
A_\bullet \rto \dto_j& M_{\bullet} \ar[d]^q \\
B_{\bullet} \ar[r] \ar@{-->}[ru] & N_{\bullet}
}
\]
can be solved in such a way that both triangles commute.
\end{prop}

We'll prove this below, but first let us record a lemma. The proof
is a diagram chase. If $M_\bullet$ is a chain complex,
let $Z_n M \subseteq M_n$ denote the cycles. We set $Z_{-1}M = 0$;
thus condition (2) in the following says $M_0 \to N_0$
is onto.

\begin{lemma}\label{lifting-lemma-1} Let  $f: M_{\bullet} \longrightarrow N_{\bullet}$ be a morphism of chain complexes. Then the following
statements are equivalent:
\begin{enumerate}

\item $H_{\ast}f$ is an isomorphism and $f: M_{n} \longrightarrow N_{n}$
is a surjection for $n > 0$.

\item The induced map
\[
M_{n} \longrightarrow Z_{n-1}M \times_{Z_{n-1}N} N_{n}
\]
is a surjection for  $n \geq 0$.
\end{enumerate}
Under either condition, the induced map $Z_{n}M \longrightarrow Z_{n}N$ 
is a surjection.
\end{lemma}  

\textsc{Proof of Proposition \ref{half-cm5}}. We would like to construct the needed
morphisms $g: B_n \to M_n$ by induction on $n$.  For the inductive
step we assume we have $g:B_k \to M_k$ for $k < n$ and that 
$g$ is a chain map as far as it is defined. Then we see that we need to
solve a lifting problem
\[
\xymatrix{
A_n \rto \dto_j& M_n \ar[d] \\
B_n \ar[r] \ar@{-->}[ru] & Z_{n-1}M \times_{Z_{n-1}N} N_{n}.
}
\]
Now apply the previous lemma and the fact that $A_n \to B_n$
is isomorphic to a morphism of the form $A_n \to A_n \oplus P$,
where $P$ is a projective.

\subsection{Model Categories}

The axioms for model categories are obtained by generalizing the
example of chain complexes. We note that in $\Ch_\ast R$ --
and in particular in Proposition \ref{half-cm5} -- we identified
three distinguished classes of maps: the homology isomorphisms,
the morphisms $M_\bullet \to N_\bullet$ which were surjections
in positive degrees, and the morphisms $A_\bullet \to B_\bullet$
which were injective with projective cokernel in all degrees. These
will become, respectively, the {\it weak equivalences}, the {\it fibrations},
and the {\it cofibrations} for $\Ch_\ast R$. As a matter of
nomenclature, we say that a morphism $f:A \to B$ is a {\it retract} of
a morphism $g:X \to Y$ if there is a commutative diagram
\[
\xymatrix{
A \dto_f \rto \ar@/^1pc/[rr]^= & X \dto^g \rto & A \dto^f\\
B \rto \ar@/_1pc/[rr]_= & Y \rto & B.
}
\]

\begin{definition}\label{model-category} A {\it model category}  is a 
category $\mathcal{C}$ with three classes of maps -- weak equivalences,  fibrations, and  cofibrations -- subject to the following axioms.
An {\it acyclic} cofibration is a cofibration which is a weak equivalence.
There is a corresponding notion of acyclic fibration.
\begin{enumerate}
\item[M1.] The category $\mathcal{C}$ is closed under limits and colimits;
\item[M2.] The three distinguished classes of maps are closed under retracts;
\item[M3.](2 out of 3) Given  $X \stackrel{f}{\longrightarrow} Y \stackrel{g}{\longrightarrow} Z$ so that any two of $f$, $g$, or $gf$ is a weak equivalence, then so is the third.
\item[M4.] (Lifting Axiom) Every lifting problem
\[
\xymatrix{
A \ar[r] \ar[d]_-{j} & X \ar[d]^-{q} \\
B \ar[r] \ar@{-->}[ru] & Y
}
\]
where $j$ is a cofibration and $q$ is a fibration has a solution so that
both diagrams commute if one of $j$ or $q$ is a weak equivalence.
\item[M5.] (Factorization) Any $f: X \longrightarrow Y$ can be factored two ways:
\begin{enumerate}
\item[(i)] $X \stackrel{i}{\longrightarrow} Z \stackrel{q}{\longrightarrow}Y$, where $i$ is a cofibration and $q$ is a weak equivalence and a fibration.
\item[(ii)] $X \stackrel{j}{\longrightarrow} Z \stackrel{p}{\longrightarrow}Y$, where $j$ is a weak equivalence and a cofibration and $p$ is a fibration.
\end{enumerate}
\end{enumerate}
\end{definition}

\begin{remarks}\label{model-cat-remarks} We immediately make some of
the standard comments about model categories.

\begin{enumerate}
\item Note that CM4 and CM5 are really two axioms each. Also note
that the axioms are completely symmetric in cofibrations and fibrations;
thus the opposite category of $\calC$ automatically inherits a model category
structure.

\item The three classes of maps are not independent.  For example,
if $\mathcal{C}$ is a model category, then the cofibrations are exactly
the morphisms with the {\it left lifting property} (LLP) with respect
to acyclic fibrations; that is, $j: A \to B$ is a cofibraton if and only if
for every acyclic fibration $q:X \to Y$ and every lifting problem
\[
\xymatrix{
A \ar[r] \ar[d]_-{j} & X \ar[d]^-{f} \\
B \ar[r] \ar@{-->}[ru] & Y
}
\]
there is a solution so that both triangles commute. In an analogous fashion,
a morphism is a fibration if and only if it has the RLP with respect to
acyclic cofibrations.
 
\item Let $X$ be an object in the model category $\calC$.
Then $X$ is {\it cofibrant} if
the unique morphism from the initial object to $X$ is a cofibration.
More generally, a {\it cofibrant replacement} or {\it cofibrant model}
for $X$ is a weak equivalence $Z \to X$
with $Z$ cofibrant. Such replacements always exist, by the factorization axiom,
and we will discuss how unique such models are below when we
talk about homotopies. (See Corollary \ref{uniqueness-of-replacements}.)

To be concrete, if $M \in \Mod_{R} \subseteq \Ch_{\ast}R$ is an
$R$-module regarded as a chain complex concentrated in
degree $0$, then a cofibrant replacement for $M$ is  simply a projective
resolution.

There are corresponding notions of {\it fibrant object} and {\it fibrant
replacement}; for example, $X$ is fibrant if the unique morphism
from $X$ to the terminal object is a fibration.

\item In his original work on the subject, Quillen required a weakened
version of M1; specifically, he only required finite limits and colimits.
He had an example, as well: chain complexes of finitely generated
$R$-modules.

\item The class of cofibrations is closed under retracts and various
colimits: coproducts,
cobase change, and sequential colimits. For example, to be closed
under cobase change means that if $i:A \to B$ is a cofibration and if
we are given any push-out diagram
$$
\xymatrix{
A \rto \dto_i & X\dto^j\\
B \rto & Y\\
}
$$
then $j$ is also a cofibration. Acyclic cofibrations are also closed
under these operations, and there are similar closure properties
for fibrations and acyclic fibrations -- now using limits in place of colimits.
\end{enumerate}
\end{remarks}

\subsection{Chain complexes form a model category}

The first example is chain complexes. For the record, we have:

\begin{theorem}\label{chains-a-model}The category $\Ch_{\ast}R$  has the structure of a model category with a morphism 
$f:M_\bullet \to N_\bullet$
\begin{enumerate}

\item a weak equivalence if $H_\ast f$ is an isomorphism;

\item a fibration if $M_n \to N_n$ is surjective for $n \geq 1$; and,

\item  a cofibration if and only if for $n \geq 0$, the map $M_n \to N_n$
is an injection with projective cokernel.
\end{enumerate}
\end{theorem}

We're going to prove this, as it gives us a chance to introduce some
notation and terminology. First note that axioms M1, M2, and M3 are
all completely straightforward. Also note that the ``acyclic fibration''
half of the lifting axiom M4 is exactly Proposition \ref{half-cm5}. The next step is
to prove the ``acyclic cofibration -- fibration'' half of the factorization
axiom M5.

To do this, let $D(n)$, $n \geq 1$ denote the chain complex with
$D(n)_k = 0$ for $k \ne n, n-1$ and
\begin{equation}\label{dn}
\partial = 1_R: D(n)_n = R \longr R = D(n)_{n-1}.
\end{equation}
Then there is a natural isomorphism
$\Ch_\ast R(D(n),N_\bullet) \cong
N_n$ given by evaluation at the generator of $D(n)_n$. Therefore,
$q:Q_\bullet \to N_\bullet$ is a fibration if and only if $q$ has the
the right lifting property with respect to the morphisms $0 \to D(n)$,
$n > 0$. 

If $N_\bullet$ is any chain complex, define a new chain complex
$P(N_\bullet)$ and an evaluation  morphim $\epsilon: P(N_\bullet)
\to N_\bullet$ by the equation
$$
P(N_\bullet) = \oplus_{n > 0} \oplus_{x \in N_n} D(n) \longr
N_\bullet
$$
This map is evidently a fibration. Better, if $M_\bullet \to N_\bullet$
is any morphism of chain complexes, we can extend this to a morphism
$$
\xymatrix{
M_\bullet \rto^-j & M_\bullet \oplus P(N_\bullet) \rto^-q & N_\bullet
}
$$
with $q$ a fibration and $j$ an acyclic cofibration with the left lifting
property with respect to all fibrations. This completes this half of
M5 and also allows us to prove the ``acyclic cofibration'' half of
M4.

Indeed, if $i: M_\bullet \to N_\bullet$ is any acyclic cofibration, the
factorization just completed yields a lifting problem
$$
\xymatrix{
M_\bullet \rto^-j \dto_i & M_\bullet \oplus P(N_\bullet) \dto^q\\
N_\bullet \rto_{=} \ar@{-->}[ur] & N_\bullet.
}
$$
Since $i$ and $j$ are weak equivalences, so is the morphism $q$; hence,
Proposition \ref{half-cm5} supplies a solution to this lifting problem
and any solution displays $i$ as a retract of $j$. Since $j$ has the
left lifting property with respect to all fibrations, so does $i$. Thus
lifting axiom M4 is verified.

This leaves the ``cofibration--acyclic fibration'' half of the factorization axion M5. We can
prove this with an induction argument. Fix a morphism of chain
complexes $f:M_\bullet \to N_\bullet$ and make the following induction
hypothesis at $n \geq 0$:

There are $R$-modules $Q_k$, $0 \leq k \leq  n-1$, and morphisms
of $R$-modules $i:M_k \to Q_k$, $p:Q_k \to N_k$ and $\partial:
Q_k \to Q_{k-1}$ so that
\begin{enumerate}
\item $pi = f: M_k \to N_k$;

\item $\partial^2 = 0$ and $i$ and $p$ are chain maps as far
as they are defined;

\item $i$ is an injection with projective cokernel and 
$$
Q_k \to Z_{k-1}Q \times_{Z_{k-1}N} N_k
$$
is a surjection for $0 \leq k \leq n-1$.
\end{enumerate}

Note that the case $n=0$ is true and easy, using the fact that all our
chain complexes
are zero in negative degrees. If we can complete the inductive step,
Lemma \ref{lifting-lemma-1} will complete M5 and the proof.
But completing the inductive step uses the same idea we used
at every other stage: that there are enough projectives. Indeed, the
maps $i$ and $q$ give a morphism
$$
f':M_n \longr Z_{n-1}Q \times_{Z_{n-1}N} N_n
$$
factoring $f$. Now choose a projective $P$ so  that there
is factoring of $f'$
$$
M_n \longr M_n \oplus P \longr Z_{n-1}Q \times_{Z_{n-1}N} N_n
$$
with the second morphism a surjection. Then set $Q_n = M_n \oplus
P$. 

As a final remark, note that with some care we could have made the
factorizations {\it natural}. See Theorem \ref{small-object}.

\subsection{Other model categories}

Here are some standard examples of model categories.

\begin{example}\label{cochain} The category $\Ch^{\ast}R$ of non-negatively graded
{\it cochain} complexes of $R$-modules supports a model category
structure which features the injective $R$-modules. The weak
equivalences are $H^{\ast}f$ isomorphisms, the cofibrations
are injections in positive degrees, and the fibrations are surjections with
an injective kernel.  Then the  fibrant replacement of a constant
cochain complex on an object $M$ is exactly an injective resolution.
The argument proceeds exactly as in the previous section.

Less obvious is that there is a model category structure
on $\ZZ$-graded chain
complexes. The argument we used for non-negatively graded chain
complexes used induction arguments which do not apply in this case.
See \cite{Model} \S 2.3. The weak equivalences are the homology
isomorphisms, the fibrations the surjections, but cofibrations
are slighly harder to understand: not every chain complex of
projectives will be cofibrant. Note that this model category structure
emphasizes projectives; there are other model category structures which
focus on injectives.
\end{example}

\begin{example}\label{formal} Suppose $\calC$ is a model category. Then
there are some formal constructions we can make to create new model
categories. For example, let $A \in \calC$ be a fixed object and
let $\calC/A$ be the over category of $A$. Thus, the objects
in $\calC/A$ are arrows $X \to A$ and the morphisms are commutative
triangles
$$
\xymatrix@C=12pt@R=15pt{
X \ar[rr]^f \ar[dr] && Y \ar[dl]\\
&A.
}
$$
The morphism $f$ is a weak equivalence, fibration, or cofibration in
$\calC/A$ if it is so in $\calC$. While this is easy, note that the fibrant
objects have changed: $q:X \to A$ in $\calC/A$ is fibrant if and only
if $q$ is a fibration in $\calC$.
\end{example}

\begin{example}\label{topology} Another basic example of a model
category is the category of topological spaces, $\mathbf{Top}$.  
A continuous map $f: X \longrightarrow Y$ is a weak equivalence if
$f_\ast:\pi_{k}(X,v) \longrightarrow \pi_{k}(Y, f(v))$ is an isomorphism
for $k \geq 0$ and all basepoints $v \in X$.\footnote{This is the origin
of the nomenclature ``weak equivalence''.} The fibrations
are Serre fibrations and Example \ref{example-small} below
implies that CW complexes are cofibrant.

The category of all topological spaces has the disadvantage that it
is not Cartesian closed: there is no exponential (or mapping space)
functor right adjoint to the product. For this reason, we often restrict to
the full subcategory $\cgh$ of compactly generated weak Hausdorff spaces. This will yield the same homotopy theory. (See the comments
after Theorem \ref{siset-and-top}.) The good
categorical properties of compactly generated Hausdorff spaces
are highlighted in the paper \cite{Steenrod}. Unforunately, the best
reference for the advantages of compactly generated weak Hausdorff spaces
remains \cite{Lewis}.
\end{example}

\subsection{Simplicial Sets}

One of the fundamental observations of non-abelian homological algebra
is that the notion of a resolution of an object by a chain complex must
be replaced by the notion of a resolution of an object by a simplicial
object. In this section, we define simplicial objects and talk about
the model category structure on simplicial sets.

For many years, the standard reference on simplicial sets has been 
\cite{May},
and this remains a good source. See also \cite{Curtis} and \cite{GZ}.
All of these references
predate model structures, so for many years we got along quite
well with \cite{BK}, Chapter VIII. The reference \cite{GJ} was
written explicitly to work model categories into the mix.

\begin{definition}\label{ordinal-numbers}The {\it ordinal number category}
$\ddelta$ has objects the ordered sets
$[n] = \{0,1, \ldots, n\}$ and  morphisms the (weakly) order preserving maps $\phi: [n] \longrightarrow [m]$. In particular we have the maps 
\begin{eqnarray*}
d^{i}: [n-1] & \longrightarrow & [n], \qquad 0 \leq i \leq n \\
\{0,1,\ldots, i-1, i, i+1, \ldots, n-1\} & \longmapsto & \{0,1,\ldots, i-1, i+1, \ldots, n\}
\end{eqnarray*}
which skip $i$ and the maps
\begin{eqnarray*}
s^{j}: [n+1] & \longrightarrow & [n] \mbox{ for } 0 \leq j \leq n \\
\{0,1,\ldots, n+1\} & \longmapsto & \{0,1,\ldots, j,j,\ldots, n\}
\end{eqnarray*}
which double up $j$. It is an exercise to show that all the morphisms
in $\mathbf{\Delta}$ are compositions of these two types of morphisms.
We also note that there are relations among them; for example,
$$
d^id^j = d^jd^{i-1},\qquad i > j.
$$
A full list of relations can be deduced from Lemma \ref{simp-relations}.
\end{definition}

\begin{definition} A {\it simplicial object} in a category $\mathcal{C}$ is a contravariant functor from $\ddelta$ to $\calC$:
\[
X:\mathbf{\Delta}^{\mathrm{op}} \longrightarrow \mathcal{C}.
\]
A morphism of simplicial objects is a natural transformation. 
The category of simplicial objects in $\mathcal{C}$ will be denoted
$s\mathcal{C}$. As a matter of notation, we
will write $X_n$ for $X([n])$,
$$
d_i = X(d^i):X_n \longr X_{n-1}\quad\hbox{and}\quad
s_i = X(s^i):X_{n-1} \longr X_n.
$$
These are respectively, the {\it face} and {\it degeneracy}
maps. More generally, if $\phi$ is morphism in $\ddelta$, we will
write $\phi^\ast$ for $X(\phi)$.

If $\calC = \mathbf{Sets}$ is the category of sets,
we get the basic category $\sisets$ of simplicial sets.
\end{definition}

The following is left as an exercise.

\begin{lemma}\label{simp-relations} For any simplicial object $X\in s\calC$, the face and degeneracy maps satisfy the following identities:
\[
\begin{array}{cccl}
d_{i}d_{j} & = & d_{j-1}d_{i} & i < j \\
d_{i}s_{j} & = & s_{j-1}d_{i} & i < j \\
           & = & 1            & i=j, j+1 \\
           & = & s_{j}d_{i-1} & i > j+1 \\
s_{i}s_{j} & = & s_{j+1}s_{i} & i \leq j.
\end{array}
\]
\end{lemma}

Here are some examples of simplicial sets.

\begin{example}[The singular set]\label{singular-set} For $n \geq 0$, let $\sigma_{n}$
be the standard topological $n$-simplex; that is, the convex hull of the
standard basis vectors in $\mathbb{R}^{n+1}$.  If $\phi:[n]
\longrightarrow [m]$ is a morphism in $\mathbf{\Delta}$, then we have a map
$$
\phi_\ast: \sigma_{n} \longrightarrow \sigma_{m}
$$
dictated by $\phi$ on the vertices in their usual ordering and extended
linearly to the rest of the simplex. Then we may define a functor
$$
S(-):\mathbf{Top} \to \sisets
$$
by setting $S_n(X) = \mathbf{Top}(\sigma_n,X)$ and $\phi^\ast = \mathbf{Top}(\phi_\ast,X)$. This is the {\it singular set} of $X$.

This object is quite familiar: if we apply the free abelian group functor
level-wise to $S(X)$ and take the alternating sum of the face maps,
we obtain the singular chain complex of first-year algebraic topology.
\end{example}

\begin{definition}\label{simplicial-nomenclature} If $X$ is a simplicial
set, we call $X_n$  the $n$-{\it simplices} and the subset
$$
\displaystyle \bigcup_{i=0}^{n-1} s_{i}X_{n-1} \subseteq X_{n}
$$
the {\it degenerate} simplices. Note that the degenerate simplices
could equally be written as
$$
\displaystyle \bigcup_{\phi:[n] \to [k]} \phi^\ast X_{k} \subseteq X_{n}
$$
where $\phi$ runs over the non-identity surjections out of $[n]$ in
$\ddelta$.
\end{definition}

\begin{example}\label{from-simplicial-complex} Let $K$ be an ordered
simplicial complex; that is, a simplicial complex with an ordering on
the set $V$ of vertices $K$. Define a simplicial set $K_\bullet$ with
$K_n$ the set of collections
$$
v_0 \leq v_1 \leq v_2 \leq \ldots \leq v_n
$$
of vertices subject to the requirement that the $v_i$, after eliminating
repetitions, are the vertices of a simplex of $K$. Thus, for example,
if $K$ is the boundary of the standard $2$-simplex with vertices
$e_i$, then $e_0 \leq e_1 \leq e_1$ is in $K_2$, but
$e_0 \leq e_1 \leq e_2$ is not. The face and
degeneracy operators are defined by doubling or deleting one of the
$v_i$.
\end{example}

\begin{example}[The standard simplices]\label{standard-simplices} The
functor from simplicial sets to sets sending $X$ to $X_n$ is 
representable, by the Yoneda Lemma. Indeed, define the standard
$n$-simplex $\Delta^n$ of $\sisets$ by
$$
\Delta^n = \ddelta(-,[n]):\Delta^{op} \to \mathbf{Sets}.
$$
Then the Yoneda Lemma supplies the isomorphism
$\sisets(\Delta^n,X) \cong X_n$. The morphisms in $\Delta$
yield morphisms $\Delta^m \to \Delta^n$.

We also have the boundary of the standard simplex
$$
\partial\Delta^n = \bigcup_{0 \leq i \leq n} d^i\Delta^{n-1} \subseteq
\Delta^n
$$
and the {\it horns}
$$
\Delta^n_k = \bigcup_{i\ne k} d^i\Delta^{n-1} \subseteq
\Delta^n.
$$
All of these simplicial sets can also be realized by the construction
of Example \ref{from-simplicial-complex} using some appropriate
subcomplex of the standard topological $n$-simplex $\sigma_n$.
\end{example}

\def\cI{{{\mathcal{I}}}}

\begin{example}[The nerve of a small category]\label{nerve}
Let $\mathcal{I}$ be a small category; that is, a category with a set
of objects. Then we define a simplicial set $B\mathcal{I}$ --
the {\it nerve}  of the category -- with the
$n$-simplices all strings of composable arrows in $\cI$:
$$
B\mathcal{I}_{n} = \{x_{0} \rightarrow x_{1} \rightarrow \cdots \rightarrow x_{n}\}.$$
If $\phi:[m] \to [n]$ is a morphism in $\ddelta$, then the 
induced function  $\phi^\ast:B\cI_n \to B\cI_m$ is given by
$$
\phi^\ast(\{x_{0} \rightarrow x_{1} \rightarrow \cdots \rightarrow x_{n}\})
= \{x_{\phi(0)} \rightarrow x_{\phi(1)} \rightarrow \cdots
\rightarrow x_{\phi(m)}\}
$$
where we use composition, deletion, or insertion of identities as
necessary. For example, the three face maps from $B\mathcal{I}_{2}$ to $B\mathcal{I}_{1}$ are given by
\[
\begin{array}{ccc}
\{x_{0} \rightarrow x_{1} \rightarrow x_{2}\} & 
\begin{array}{c} \overset{d_0}{\longmapsto} \\ \overset{d_1}{\longmapsto }\\ \overset{d_2}{\longmapsto} \end{array} &
\begin{array}{c}
\{x_{1} \rightarrow x_{2}\} \\
\{x_{0} \rightarrow x_{2} \}\\
\{x_{0} \rightarrow x_{1}\}.
\end{array}
\end{array}
\]
To be concrete, if $\cI$ is the category
$$
0 \to 1 \to \cdots \to n,
$$
then $B\cI \cong \Delta^n$.
\end{example}  

\begin{example}\label{the-1-simplex} To make your confusion completely
specific\footnote{This joke is stolen from Steve Wilson.},
consider $\Delta^1$ as the nerve of $0 \to 1$. 
Then we have
$$
\mathbf{\Delta}^{1}_{0} = \{ 0, 1\}
$$
and
$$
\mathbf{\Delta}^{1}_{1} =  \{ 0 \to 0, 0\to 1, 1\to 1\}
$$
and
$$
\mathbf{\Delta}^{1}_{2} = \{0\to 0\to 0, 0\to 0\to 1, 0\to 1\to 1\}.
$$
Thus the only non-degenerate simplex
in $\mathbf{\Delta}^{1}_{1}$ is $0 \to 1$ and all of the simplices in $\mathbf{\Delta}^{1}_{2}$ are degenerate.
\end{example}

\begin{example}[The product of simplicial sets] Unlike the product
of simplicial complexes, the product of simplicial sets is easy to describe.
If $X$ and $Y$ are simplicial sets, then $(X \times Y)_{n} = 
X_{n} \times Y_{n}$ and $\phi$ is a morphism in $\ddelta$, then
$$
\phi^\ast_{(X \times Y)} = \phi^\ast_X \times \phi^\ast_Y.
$$
\end{example}

\begin{definition}\label{geometric-realization} For a simplicial set $X$
define a topological space $|X|$ -- the {\it geometric realization} of
$X$ -- by taking a coequalizer in the category $\cgh$ of compactly
generated weak Hausdorff spaces (see Example \ref{topology}):
\[
\xymatrix{
\coprod_{\phi: [n] \rightarrow [m]} X_{m} \times \sigma_{n}  \ar@<2pt> [r] \ar@<-2pt> [r] & \coprod_{n} X_{n} \times \sigma_{n}  \ar [r] & |X| 
}
\]
Here $\sigma_n$ is the topological $n$-simplex and
$X_n \times \sigma_n = \coprod_{X_n} \sigma_n$. The
parallel arrows are the two arrows that are induced by evaluating
on $X$ or on $\sigma_{(-)}$ respectively. 
\end{definition}

The following result tabulates some facts about
the geometric realization functor. The singular set functor
was defined in Example \ref{singular-set}. In the second
point we are using the product in $\cgh$ -- the ``Kelly
product''. See \cite{Steenrod}.

\begin{proposition}\label{props-realization} The geometric
realization functor has the following properties:
\begin{enumerate}
\item The functor $|-|: \sisets \to \cgh$ is left
adjoint to the singular set functor $S(-)$.

\item If $X$ and $Y$ are simplicial sets, then the natural morphism
$|X \times Y| \to |X| \times |Y|$ in $\cgh$ is a homeomorphism.

\item If $\Delta^n$ is the simplicial $n$-simplex, then
$|\Delta^{n}| \cong \sigma_{n}$.
\end{enumerate}
\end{proposition}

\begin{remark}\label{drawing} We use geometric realization
to depict simplicial sets. Only the non-degenerate simplices
are needed for a picture of $|X|$, and the face maps tell
us how to glue them together.
For example, for $\mathbf{\Delta}^{1} \times \mathbf{\Delta}^{1}$ we have the following picture of the geometric realization:
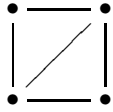
\[ 
\xymatrix{
\bullet \ar@{-} [r] \ar@{-} [d] & \bullet \ar@{-} [d] \\
\bullet \ar@{-} [ru] \ar@{-} [r]  & \bullet
}
\]
There are four non-degenerate $0$-simplices, five non-degenerate
$1$-simplices and two non-degenerate $2$-simplices glued as
shown. (There are also orientations on the $1$- and $2$-simplices,
which are not shown.)
Note that while $\Delta^1$ has no non-degenerate $2$-simplices,
$\Delta^1 \times \Delta^1$ indeed has two non-degenerate $2$-simplices.
It is an extremely useful exercise to identify these elements of
$(\Delta^1 \times \Delta^1)_2$.
\end{remark}

We now come to the model category structure on simplicial sets.
It is still one of the deeper results in the theory. I know of no easy
proof, nor of one that doesn't essentially also prove the Quillen
equivalence of simplicial sets and topological spaces at
the same time. (See Theorem \ref{siset-and-top} below.) The arguments
are spelled out in \cite{GJ} and \cite{Model} but these
are variations on Quillen's original argument of \cite{HA}.

\begin{theorem}\label{model-simplicial} The category $\sisets$
has the structure of a model category with a morphism
$f:X \to Y$
\begin{enumerate}
\item a weak equivalence if  $|f|:|X| \to |Y|$ is a weak equivalence
of topological spaces;

\item a cofibration if $f_n:X_n \to Y_n$ is a monomorphism for all $n$; and

\item a fibration if $f$ has the left lifting property with respect to the
inclusions of the horns
$$
\Delta^n_k \longr \Delta^n, \qquad n \geq 1,\ 0 \leq k \leq n.
$$
\end{enumerate}
\end{theorem}

The lifting condition in the definition of fibration for simplicial
sets is known as the {\it Kan extension condition} and goes back
to Kan's original work on the subject in \cite{comb}. It arose there
as a condition that made it possible to define and analyze the homotopy
groups of a simplicial set without referring to topological spaces. (Compare
Example \ref{moore-complex}.)
However, the lifting condition implies that a fibrant simplicial set has very many simplices,
which makes them hard to draw: for example, $\Delta^n$ is {\it not}
fibrant if $n > 0$. Or, for another example, if $\cI$ is a small category,
its nerve  $B\cI$ is fibrant if and only if $\cI$ is a groupoid.

For some good news, however, we note that every simplicial group, regarded as
a simplicial set, is fibrant; indeed, every surjective morphism
of simplicial groups is a fibration. To make a definitive statement
along these lines, let us make a definition. If $X$ is a simplicial
set, define $\pi_0X$ by the coequalizer diagram of sets
\begin{equation}\label{pi0}
\xymatrix{
X_1 \ar@<.5ex>[r]^{d_1} \ar@<-.5ex>[r]_{d_0} &X_0 \rto &\pi_0X.
}
\end{equation}
It is an exercise to show that $\pi_0X \cong \pi_0|X|$. The
following is also left to the reader.

\begin{prop}\label{fib-sim-grps}Let $f:X \to Y$ be a morphism
of simplicial groups. Then $f$, regarded as a morphism of simplicial
sets, is a fibration if and only if
$$
X \longr \pi_0X \times_{\pi_0Y} Y
$$
is a surjection.
\end{prop}

\section{Quillen Functors and Derived Functors}

We have defined the objects of interest -- model categories -- and
we will 
now define and study the morphisms between them, which are Quillen
functors. Then we will introduce the notion of the homotopy
category of a model category and discuss the fact that Quillen functors
descend to adjoint functors between the homotopy categories. This is
a little backwards -- traditionally one defined the homotopy category 
almost immediately after defining the model category -- but the lesson of
the last thirty or so years is that there's more information in the model
category than in the homotopy category and we should descend later
rather than sooner.

\subsection{Quillen Functors and Quillen Equivalence}

Here are our morphisms between model categories.

\begin{definition}
Let $\mathcal{C}$ and $\mathcal{D}$ be two model categories. Then
a {\it Quillen functor} from $\calC$ to $\mathcal{D}$ is an adjoint pair of functors
\[
\xymatrix{
\mathcal{C} \ar@<2pt> [r]^-{F} & \mathcal{D} \ar@<2pt>[l]^-{G}
}
\]
with $F$ the left adjoint so that
\begin{enumerate}
\item the functor $F$ preserves cofibrations and weak equivalences between cofibrant objects, and 
\item the functor $G$ preserves fibrations and weak equivalences between fibrant objects.
\end{enumerate}
A Quillen functor is a {\it Quillen equivalence} if for all cofibrant objects
$X$ in $\calC$ and all fibrant objects $Y$ in $\mathcal{D}$, a morphism
\[
X \longr GY
\]
is a weak equivalence in $\mathcal{C}$ if and only if the adjoint morphism
\[
FX \longr Y
\]
is a weak equivalence in $\calD$.
\end{definition}

Here and elsewhere we will write an adjoint pair of functors with the left
adjoint on top, and we regard this as a ``function'' from $\calC$ to $\calD$.

\def\res{{{\mathrm{res}}}}

\begin{example}\label{tensor-as-quillen}
Let $f:R \longrightarrow S$ be a homomorphism of commutative
rings and let $\res_f$ denote
the restriction of scalars functor from $S$-modules to $R$-modules.
\[
\xymatrix{
S\otimes_{R} -: \Ch_{\ast}R \ar@<2pt>[r] & \Ch_{\ast}S:\res_f \ar@<2pt> [l]}
\]
yields a Quillen functor. Note that the functor $\res_f$ preserves
all weak equivalences. This is a Quillen equivalence only if
$R=S$.
\end{example}

The following liberating result, which first appeared in \cite{HA},
has its roots in the work of Kan and is a modern formulation of
Kan's observation that the homotopy theory of topological spaces
and the homotopy theory of simplicial sets are equivalent.

\begin{theorem}\label{siset-and-top}
The geometric realization functor and the singular set functor
give a Quillen equivalence
\[
\xymatrix{
|-|: \sisets \ar@<2pt> [r] & \cgh:S(-). \ar@<2pt> [l]}
\]
\end{theorem}

One might quibble that one would really want the homotopy theory
of all topological spaces, not only compactly generated weak
Hausdorff spaces. However, a result we leave the reader to formulate yields
an appropriate Quillen equivalence between the two. But we should also
append the remark that we are using the {\it Serre} model category
structure on spaces, as in Example \ref{topology}. There is
another model category structure on spaces wherein the weak equivalences
are actually homotopy equivalences. See \cite{Strom}. This is another
story all together, and we leave it for another day.
 
\subsection{Homotopies and the homotopy category}

So far we have been writing about model categories as if they encoded
a homotopy theory. This is indeed the case, and in this section, we make
this precise.

\begin{definition}\label{cylinder}
Let $A \in \mathcal{C}$, where $\mathcal{C}$ is a model category.  A \emph{cylinder object} for $A$  is a factoring
\[
\xymatrix{
A \amalg A \ar [r]^-{i} \ar [dr]_-{\nabla} & C(A) \ar[d]^-{q} \\
& A
}
\]
where $i$ is a cofibration, $q$ is a weak equivalence, and $\nabla$ is
the fold map. We will refer to $C(A)$ as the cylinder object, leaving the
rest of the diagram implicit. 
\end{definition}

\begin{example}\label{topology-cylinder} If $A$ is a CW complex,
then we may choose $C(A) = A \times [0,1]$. But note that if
$A$ is not cofibrant
as a topological space, than the inclusion
$$
A \amalg A \longr A \times [0,1]
$$
may not be a cofibration. We also remark that at this point that this is
an example of a {\it natural} cylinder object: this is desirable, and may
even be required at certain points. Compare Theorem \ref{quillen-criterion}.
\end{example}

\begin{example}\label{chain-cylinder} We may define a {\it natural} cylinder
object for chain complexes as follows. Let $P_{\bullet} \in \Ch_{\ast}R$ be cofibrant, and set
$$
C(P_{\bullet})_{n} = P_{n} \oplus P_{n-1} \oplus P_{n}
$$
with boundary 
$$
\partial (x,a,y) = (\partial x + a, -\partial a, \partial y + a).
$$
In all the examples of this monograph, we will have natural
cylinder objects, and we could have included this in our
definitions. 
\end{example}

The following definition is motivated by the example of topological
spaces.

\begin{definition}\label{left-homotopy} Let $f,g:A \longrightarrow X$ be
two morphisms in a model category $\calC$. 
A \emph{left homotopy} from $f$ to $g$ is a diagram in $\calC$
\[
\xymatrix{
A \sqcup A \ar[r]^i \ar[dr]_-{f \sqcup g} & C(A) \ar [d]^-{H} \\
& X.
}
\]
where $C(A)$ is a cylinder object for $A$. We will denote the 
homotopy by $H$ and drop the ``left'' when it can be understood
from the context.
\end{definition}

\begin{remark}\label{right-homotopy} We leave it to the reader to
formulate the appropriate notions of a {\it path object} and {\it
right homotopy}. These ideas are completely dual to those just
presented.
\end{remark}

\begin{example}\label{chain-homotopy}Using the cylinder object
for chain complexes we wrote down in Example \ref{chain-cylinder}
we see that in $\Ch_{\ast}R$, the homotopies correspond to the chain homotopies.
\end{example}

\begin{example}\label{constant-homotopy} Here is a formal example.
Given $f:A \to X$, a {\it constant homotopy} is given by the diagram
\[
\xymatrix{
A \sqcup A \ar [r]^i \ar [d]_-{f \sqcup f} & C(A) \ar [d]^q \\
X & A \ar [l]^-{f}.
}
\]
If $\varphi:X \to Y$ is a morphism in $\calC$ and $f,g:A \to X$ are two
morphisms so that $\varphi f = \varphi g$, then a homotopy from $f$ to
$g$ {\it over Y} is a homotopy $H:C(A) \to X$ so that $\varphi H$
is the constant homotopy. There is also the notion of a right homotopy
under an object.
\end{example}

The following result now follows from the model category axioms.
There is a corresponding result for right homotopies.

\begin{lemma}\label{uniqueness-of-lifting} 
Suppose there is a lifting problem in a model category $\calC$
\[
\xymatrix{
A \ar [r] \ar[d]_-{j} & X \ar[d]^-{q} \\
B \ar@{-->} [ur]^-{f} \ar [r] & Y
}
\]
where $j$ is a cofibration, $q$ is a fibration, and $q$ is a
weak equivalence. Then $f$ exists and $f$ is unique up to left homotopy
under $A$ and over $Y$.  In other words, any two solutions to the lifting
problem are homotopic via a homotopy which reduces to a constant homotopy when restricted to $A$ or pushed forward to $Y$.

If $j$ is a weak equivalence and $A$ is cofibrant, any two lifitngs are left
homotopic under $A$ and over $Y$.
\end{lemma}

The next result now follows immediately and provides a uniqueness
result for cofibrant replacements. Note that the model category
axioms imply that for any object $A$ in a model category there
is an acyclic fibration $X \to A$ with $X$ cofibrant.

\begin{corollary}\label{uniqueness-of-replacements}Let $\calC$ be a model
category and $A$ an object in $\calC$. Suppose we are given
a weak equivalence $f_1:X_1 \to A$ and 
an acyclic fibration  $f_2:X_2 \to A$ with
$X_1$ and $X_2$ cofibrant. Then there is a weak equivalence
$X_1 \to X_2$ over $A$ and this morphism is unique up to homotopy
over $A$. Furthermore, if $f_1$ is an acyclic fibration, then 
$X_1$ and $X_2$ are homotopy equivalent over $A$.
\end{corollary}

As before we leave the reader to formulate a version of this
result for fibrant replacements.

\def\Ho{{\mathbf{Ho}}}

\begin{definition}[The homotopy category]\label{homotopy-category}
Let $\mathcal{C}$ be a model category. Then the {\it homotopy
category}  $\Ho(\calC)$ is the category obtained from $\calC$ by
inverting the weak equivalences. Thus $\Ho(\calC)$ is characterized
by the property that there is a functor $\iota:\calC \to \Ho(\calC)$ 
which takes weak equivalences to isomorphisms and if
$F:\calC \to \calD$ is any functor taking weak equivalences to isomorphisms
then there is a unique functor $F':\Ho(\calC) \to \calD$
making the following diagram commute
$$
\xymatrix{
C \rto^-\iota \dto_F & \Ho(\calC) \ar[dl]^{F'}\\
\calD\\
}
$$
We may write $\Ho(\calC)(X,Y)$ as $[X,Y]_\calC$.
\end{definition}

The morphisms in $\Ho(\calC)$ can be described as zig-zags where
the reverse arrows are weak equivalences, but this begs the question
of whether the morphism sets are actually sets. It is more practical
to give the following two descriptions of the morphisms. First,
there is an isomorphism 
\begin{equation}\label{homotopy-classes}
\Ho(\mathcal{C})(X,Y) \cong
{\mathcal{C}}(X_{c}, Y_{f})/(\sim \mathrm{homotopy}),
\end{equation}
where $X_{c} \to X$ is a cofibrant replacement and
$Y \to Y_{f}$ is fibrant replacement. The fact that the right hand side
of this equation is well-defined follows from Corollary 
\ref{uniqueness-of-replacements} and its fibrant analog.

Second, the following result, due to Dwyer and Kan \cite{hammock},
gives a simple choice-free description of the morphisms, although
it also has its set-theoretic issues.

\begin{lemma} Let $\calC$ be a model category and let $X$ and $Y$ be
objects in $\calC$.  Then $\Ho(\calC)(X,Y)$ is isomorphic to the set
of equivalence classes of diagrams in $\calC$
$$
X\stackrel{\simeq}{\longleftarrow} U \longrightarrow V \stackrel{\simeq}{\longleftarrow}Y 
$$
where the left-facing arrows are weak equivalences. The equivalence
relation among these diagrams is the smallest equivalence relation 
containing the relation created by diagrams of the form
$$
\xymatrix@R=5pt{
&\ar[dl]_-\simeq U_1 \ar[dd] \rto & V_1\ar[dd]\\
X&&&\ar[ul]_-\simeq \ar[dl]^-\simeq Y. \\
& \ar[ul]^-\simeq U_2 \rto & V_2
}
$$
\end{lemma}

As the reader may have surmised, this description of the morphisms in
the homotopy category can be obtained by taking the components of
the nerve of a suitable category. 

\def\sisphere#1{{{\Delta^{{#1}}/\partial\Delta^{{#1}}}}}
\def\sihorn#1{{{\Delta^{{#1}}/\Delta_0^{{#1}}}}}

\begin{example}[Homotopy groups of simplicial sets]\label{moore-complex} In simplicial sets, we can take $\sisphere{n}$ as
a model for the $n$-sphere. Hence, if $X$ is a pointed fibrant 
simplicial set, $\pi_n X$ can be calculated as the based homotopy classes
of maps $\sisphere{n} \to X$. Since the geometric realization of
$\sisphere{n}$ is the topological $n$-sphere, Theorem
\ref{siset-and-top} implies that these homotopy groups are naturally
isomorphic to the homotopy groups of the geometric realization
of $X$. This yields Kan's combinatorial
definition of homotopy groups. See \cite{comb}.

If $X$ is a simplicial {\it group}, then $X$ is automatically fibrant
and there is a natural isomorphism
$$
\sisets_\ast (\sisphere{n},X) \cong \bigcap_{i=0}^{n}\mathrm{Ker}\{d_{i}: X_{n} \longrightarrow X_{n-1} \}
$$
(Here the asterisk means were are taking based maps; $e \in X$ is the basepoint
of the simplicial group $X$.)
Furthermore, a morphism $\sisphere{n} \to X$ is null-homotopic
if and only if it can be factored
$$
\xymatrix{
\sisphere{n} \rto^-{d^0} & \sihorn{n+1} \rto &X.
}
$$
Hence, if we set
$$
NX_n = \bigcap_{i=1}^{n}\mathrm{Ker}\{d_{i}: X_{n} \longrightarrow X_{n-1} \}
$$
the homomorphisms $d_0:NX_{n+1} \to NX_{n}$ satisfy $d_0^2 = 0$
and we can conclude that $\pi_nX \cong H_n(NX)$. The chain
complex of not-necessarily abelian groups $(NX,d_0)$ is
called the {\it Moore complex}. If $X$ is a simplicial abelian
group, $NX$ is also called the {\it normalized chain complex}
of $X$. See \S 4.1 for more discussion.
\end{example}

\subsection{Total derived functors}

Returning to the thread of model categories as a source of resolutions,
we now use model categories to define derived functors.

\begin{definition}\label{total-left} Let $\xymatrix{F:\mathcal{C} \ar@<2pt> [r] & \mathcal{D}:G \ar@<2pt> [l]}$ be a Quillen functor.  Then $F$ has a \emph{total left derived functor} $LF$ defined as follows.  If $X \in \mathcal{C}$, let $X_c \to X$ be a weak equivalence with $X_c$ cofibrant. 
Then set $LF(X) = F(X_c)$.
\end{definition}

It immediately follows from Corollary \ref{uniqueness-of-replacements}
that $LF(X)$ is independent, up to weak equivalence, of the choice of
cofibrant replacement, and is well-defined in the homotopy category of $\calD$. 
As a side remark, we can weaken the assumption that we have a Quillen
functor. Quillen, for example, describes the total left derived functor
as a Kan extension of sorts. See \cite{HA} or \cite{GJ} \S II.7.

Here is a first example. There will be others below.

\begin{example}\label{tor-is-total} Let $f:R \to S$ be a morphism of
commutative rings and consider the Quillen functor of 
Example \ref{tensor-as-quillen}
$$
\xymatrix{S\otimes_{R} - :\Ch_{\ast}R \ar@<2pt> [r] & \Ch_{\ast}S:\res_f. \ar@<2pt>[l]
}
$$
Let $M \in \Mod_{R} \subseteq \Ch^{\ast}R$ and
 $P_{\bullet} \stackrel{\cong}{\longrightarrow} M$ be a projective resolution (i.e., a cofibrant replacement) of $M$.  Then
\[
L(S\otimes_{R}-)(M) = S\otimes_R P_{\bullet}.
\]
In particular,
$$
H_nL(S\otimes_{R}-)(M) = \mathrm{Tor}_n^R(S,M).
$$
For more general objects $N_\bullet \in \Ch_\ast R$, there is a spectral 
sequence
$$
\mathrm{Tor}_p^R(S,H_qN_\bullet) \Longrightarrow 
H_{p+q}L(S\otimes_{R}-)(N_\bullet).
$$
See \cite{HA}, \S II.6 or Example \ref{ss-simp-model}. Note that
we often write
$$
L(S \otimes_{R}-)(N_\bullet) = S \otimes_R^L N_\bullet.
$$
\end{example}

There is also a companion notion of total right derived functor. Quillen
introduced the notion of a Quillen functor (without using that name,
of course) in order to prove the following result -- which applies,
in particular, to the case of simplicial sets and topological spaces.
See Theorem \ref{siset-and-top}.

\begin{proposition}\label{quillen-homotopy}  Let  $\xymatrix{F:\mathcal{C} \ar@<2pt> [r] & \mathcal{D}:G \ar@<2pt> [l]}$ be a Quillen functor
between two model categories. Then the total derived functors induce
an adjoint pair
$$
\xymatrix{
LF:\Ho(\calC) \ar@<2pt> [r] & \Ho(\calD) :RG. \ar@<2pt> [l]
}
$$
Furthermore, this adjoint pair induces an equivalence of categories
if and only if the Quillen functor is a Quillen equivalence.
\end{proposition}

\section{Generating New Model Categories}

In this section we discuss how to promote model category structures
from one category to another. The basic idea is to take an adjoint
pair $\xymatrix{F:\mathcal{C} \ar@<2pt> [r] & \mathcal{D}:G \ar@<2pt> [l]}$
and a model category structure on $\calC$, and then force it
to be a Quillen functor. To do this requires structure and hypotheses,
which we explore first.

\subsection{Cofibrantly Generated Model Categories}

\begin{definition}\label{small} Let $\mathcal{C}$ be a category and let
$\mathcal{F} \subseteq \mathcal{C}$ be a class of maps.
Then $A \in \mathcal{C}$ is {\it small} for $\mathcal{F}$ if
whenever
$$
X_{1} \rightarrow X_{2} \rightarrow X_{3} \rightarrow \cdots
$$
is a sequence of morphisms in $\mathcal{F}$, then the natural map
\[
\mathrm{colim}\ \mathcal{C}(A,X_{n}) \longrightarrow 
\mathcal{C}(A, \mathrm{colim} X_{n})
\]
is an isomorphism.
\end{definition}

Here we are taking colimits over the natural numbers. Once could
define a notion of $\lambda$-small where $\lambda$ is any ordinal
number. This is useful for localization theory (see \cite{Local}), but
won't appear in any of our examples here. 

\begin{examples}\label{example-small} Here are some basic examples.
\begin{enumerate}
\item In the category $\cgh$ of compactly generated Hausdorff spaces,
a compact space $A$ is $J$-small, when $J$ is the class
of closed inclusions. See Proposition 2.4.2 of \cite{Model}.

\item Any bounded chain complex of finitely presented $R$-modules is $J$-small, for $J$ all morphisms of chain complexes.

\item Any finite simplicial set $X$ -- meaning $X$ has finitely many non-degenerate simplices --  is $J$-small, for $J$ the class of all morphisms
of simplicial sets. This is one of the reasons that simplicial sets are
technically pleasant.
\end{enumerate}
\end{examples}

\begin{definition}\label{cof-gen} A model category is {\it cofibrantly
generated} if there are sets of morphisms $I$ and $J$ so that
\begin{enumerate}

\item the source of every morphism in $I$ is small for the class of all
cofibrations and $q: X \longrightarrow Y$ is an acyclic fibration if and
only if $q$ has the RLP with respect to all morphisms of $I$; and

\item the source of every morphism in $J$ is small with respect to the
class of all acyclic cofibrations and $q:X \longrightarrow Y$ is a
fibration if and only if $q$ has the RLP with respect to all morphisms of
$J$. 
\end{enumerate}
\end{definition}

The set $I$ and the set $J$ generate the cofibrations and the acyclic
cofibrations repsectively.  Here, ``generate'' means the cofibrations are
the smallest class of maps that contains $I$ and is closed under coproducts,
cobase change, sequential colimits, and retracts.  

\begin{examples}\label{exam-cof-gen}All of our basic model
categories are cofibrantly generated.

1. First consider the category  $\Ch_{\ast}R$ of chain
complexes over a commutative ring $R$. The set $J$ of generating
acyclic cofibrations can be taken to be the morphisms
$$
0 \longr D(n),\qquad n \geq 1
$$
where $D(n)$ is the object so that $\Ch_\ast R(D(n),M_\bullet)
\cong M_n$. (See Equation \ref{dn}.) To specify a set of generating
cofibrations, let $S(n)$ be the chain complex with $S(n)_n = R$
and $S(n)_k = 0$ for $k \ne n$. Then the set $I$ can be chosen to consist
of the obvious inclusions
\begin{align*}
S(n-1) &\longr D(n),\qquad n \geq 1;\\
0 &\longr S(0)
\end{align*}
This follows immediately from Proposition \ref{half-cm5}.

2. The Serre model category on compactly generated weak Hausdorff
spaces is cofibrantly generated. We may choose the set $I$ of generating
cofibrations to be inclusions of boundaries of the disks
$$
S^{n-1} = \partial D^n \longrightarrow D^{n},\qquad n \geq 0
$$ 
where $S^{-1} = \emptyset$ and we may choose $J$ to be the inclusions
$$
D^{n} \longrightarrow D^{n} \times [0,1],\qquad n \geq 0.
$$
Indeed, to have the RLP with respect to the elements of $J$ is
exactly the definition of a Serre fibration.

3. The category of simplicial sets
 is cofibrantly generated.  The generating cofibrations can be
chosen to be the inclusions
$$
\partial\Delta^{n} \longr \Delta^{n}, \qquad n \geq 0
$$
where, again $\partial\Delta^0 = \emptyset$ and the generating
acyclic cofibrations are the inclusions of the horns\footnote{
Also known as anodyne extensions. See \cite{GZ}.}
$$
\Delta^{n}_{k} \longr \Delta^{n},\qquad 0 \leq k \leq n,\ n \geq 1.
$$
The horns and $\partial\Delta^n$ were defined in Example \ref{standard-simplices}.

4. But not every model category is cofibrantly generated. Perhaps the
easiest example, admittedly a bit artificial, is the induced model category
structure on the opposite category of simplicial sets.
\end{examples}

We now come to a basic fact about cofibrantly generated model
categories. The proof of this fact may be more important than the
result itself. Called Quillen's {\it small object argument}, it
first appeared in \cite{HA}, \S II.4, and we present it in some detail.

\begin{theorem}[The small object argument]\label{small-object}
If $\mathcal{C}$ is cofibrantly generated, then factorizations
of M5 can be chosen to be natural.
\end{theorem}

\begin{proof} We will do the 	``cofibration-acylic fibration'' factorization;
the other is similar. Let $f: X \longrightarrow Y$ in $\mathcal{C}$.
We want to factor $f$ in a natural way as
\[
\xymatrix{X \ar[r]^-{j} & Z \ar[r]^-{q} & Y}
\]
where $j$ is a cofibration and $q$ is an acyclic fibration.  To do this,
we produce the diagram over the object $Y$
\[
\xymatrix{
X=Z_{0} \ar [r]^-{j_{1}} \ar [rd]_-{f=q_{0}} & Z_{1} \ar [r]^-{j_{2}} \ar [d]^-{q_{1}} & Z_{2} \ar [r] \ar[dl]^-{q_{2}} & \cdots \\
& Y 
}
\]
so that there is a pushout diagram
\[
\xymatrix{
\amalg_{U} A \ar[r] \ar[d] & Z_{n} \ar[d]^{j_n} \\
\amalg_{U} B \ar[r] & Z_{n+1}
}
\]
where $U$ is the set of diagrams of the following form\[
\xymatrix{
A \ar[r] \ar[d]_-{i} & Z_{n} \ar[d]^-{q_{n}} \\
B \ar[r] & Y
}
\]
with $i \in I$.  Then $Z = \mathrm{colim}\, Z_{n}$ and $q = \mathrm{colim}\, q_{n}$.  The induced map $j:X \to Z$
is a cofibration because it is generated naturally by elements in $I$.  
Thus we need only show that $q$ is an acyclic fibration. For this, it 
is sufficient  to find a solution to every lifting problem
\[
\xymatrix{
A \ar[r]^{\phi} \ar[d]_-{i} & Z \ar[d]^-{q} \\
B \ar[r]_\psi \ar@{-->}[ru] & Y
}
\]
where $i \in I$.  Because $A$ is small with respect to the class of all
cofibrations, there is a factoring
\[
\xymatrix{
& Z_{n} \ar[d] \\
A \ar[ur]^-{\tilde{\phi}} \ar[r]^-{\phi} & Z
}
\]
for some $n$; therefore,  by the construction of $Z_{n+1}$,
we have a diagram 
\[
\xymatrix{
A \ar[r]^-{\tilde{\phi}} \ar[d] & Z_{n} \ar[d] \\
B \ar[r] \ar[dr]_-{\tilde{\psi}} & Z_{n+1} \ar[d] \\
& Z
}
\]
and $\tilde{\psi}$ solves the original lifting problem.
\end{proof}

\subsection{Promoting model category structures}

We now give a result -- again essentially due to Quillen -- for lifting
a model category structure from one category to another. Recall from
Remark \ref{model-cat-remarks}.2 that once we have specified weak
equivalences and fibrations, then cofibrations will be forced. Thus
a ``cofibration'' in the following statement means a morphism with
the left lifting property with respect to all acyclic fibrations.

\begin{theorem}\label{lifting} Let
$\xymatrix{F: \mathcal{C} \ar@<2pt>[r] & \mathcal{D}:G
\ar@<2pt>[l]}$ be an adjoint pair and suppose $\mathcal{C}$ is a
cofibrantly generated model category.  Let $I$ and $J$ be chosen sets
of generating cofibrations and acyclic cofibrations, respectively. Define
a morphism  $f: X \to Y$ in $\calD$ to be a weak equivalence or 
a fibration if $Gf$ is a weak equivalence or fibration in $\calC$.
Suppose further that
\begin{enumerate}

\item the right adjoint $G:\calD \to \calC$ commutes with sequential colimits; and

\item every cofibration in $\mathcal{D}$ with the LLP with respect to
all fibrations is a weak equivalence.
\end{enumerate}

Then $\mathcal{D}$ becomes a cofibrantly generated model category.
Furthermore the sets $\{\ Fi\ |\ i \in I\ \}$ and $\{\ Fj\ |\ j \in J\ \}$
generate the cofibrations and the acyclic cofibrations of $\calD$
respectively.
\end{theorem}

The proof, which makes repeated use of the small object argument,
has been presented many places. See \cite{GJ}, Chapter II, for example. 
See also \cite{Local} and \cite{Model}. But all these arguments
are riffs of the argument in \cite{HA}, \S II.4.

\def\dga{{{\mathbf{DGA}}}}

\begin{example}\label{something-to-prove} Some hypothesis is
needed, as it is not always possible to lift model category structures
in this way. 
If the conclusion of Theorem \ref{lifting} holds and the category $\calD$
does have the indicated model category structure, then a simple
adjointness argument shows that the left adjoint $F$ must preserve
all weak equivalences between cofibrant objects. We give an example
where this does not hold.

Let $\calC = \Ch_\ast k$ where $k$ is a field of characteristic
$2$ and let $\calD = \dga_k$ denote the category of non-negatively
graded commutative differential graded algebras. Then the forgetful
functor $\dga_k \to \Ch_\ast k$ has a left adjoint $S$ given by the
the symmetric algebra functor, with differential extended by the
Leibniz rule. Every object in $\Ch_\ast k$ is cofibrant.
If $D(n) \in \Ch_\ast$ is the chain complex with
one copy of $k$ in degrees $n$ and $n-1$ and identity boundary,
then $0 \to D(n)$ is a weak equivalence. However, $S(0) = k \to
S(D(n))$ is not a weak equivalence. Indeed if $y \in D(n)_n$ is
non-zero, then $y^2 \in S(D(n))$ is a cycle of degree $2n$
which is not a boundary.

If, on the other hand, $k$ is a field of characteristic zero, then
$\dga_k$ does inherit a model category structure in this way.
\end{example}

There is a useful criterion for checking hypothesis (2) of Theorem
\ref{lifting}, also due to Quillen in his original work \cite{HA} and
recently highlighted in \cite{BergerM}.

\begin{theorem}\label{quillen-criterion} Suppose in the category
$\mathcal{D}$ (which is not yet a model category) the following two
conditions hold:
\begin{enumerate}
\item there is a functorial fibrant replacement functor; and
\item every object has a natural path object. In other words we have a
natural diagram for all $B \in \mathcal{D}$
\[
\xymatrix@C=30pt{
& P(B) \ar[d]^q \\
B \ar[ur]^-{i} \ar[r]_-{\Delta} & B \times B.
}
\]
where $i$ is a weak equivalence and $q$ is a fibration.
\end{enumerate}
Then every cofibration in $\mathcal{D}$ with the LLP with respect to
all fibrations is a weak equivalence.
\end{theorem}

There is a standard situation where path objects automatically
exist -- namely, when $\calD$ is a simplicial category. See
Corollary \ref{path-cylinder}.

\section{Simplicial Algebras and Resolutions in Non-abelian Settings}

\subsection{The Dold-Kan theorem and simplicial resolutions}

\def\Alg{{{\mathbf{Alg}}}}

In Example \ref{something-to-prove}, we showed that the free
commutative $R$-algebra functor, extending to chain complexes,
does not generally preserve weak equivalences. As a consequence,
differential graded algebras do not provide a good setting
for creating resolutions of commutative algebras -- except in
characteristic zero. However, quite
early on, Dold \cite{dold} noticed that the free commutative
algebra functor {\it did} preserve weak equivalences between cofibrant
simplicial $R$-modules and that, therefore, simplicial algebras
can provide a suitable setting for resolutions. This was extended
and discussed further by Dold and Puppe in \cite{dold-puppe},
which, in effect, talked about total derived functors before Quillen
put a name to that concept.

In this section, we will begin with a discussion of simplicial
$R$-modules in order to set the stage for discussing other
sorts of categories of simplicial objects. The observation is that
the category of simplicial $R$-modules is equivalent to the 
category of non-negatively graded chain complexes over $R$,
so that one can work equally well with either category. Once this
discussion is in place, we'll define a notion of resolution of
a commutative algebra.

Let $s\Mod_R$ be the category of simplicial $R$-modules and $X \in
s\Mod_R$. Then we defined the {\it normalized} chain complex
of $X$ by eliminating the degenerate simplices as follows:
$$
NX_n = \frac{X_{n}}{s_{0}X_{n-1} + \cdots + s_{n-1}X_{n-1}}
$$
and  setting
$$
\partial = \sum_0^n (-1)^n d_i:NX_n \longr NX_{n-1}.
$$
It is an exercise to show that the homomorphism $\partial$ 
is well defined and that
$$
H_\ast NX \cong  H_\ast (X,\sum_0^n (-1)^n d_i).
$$ 
There is an isomorphic formulation of $NX$. (See Example
\ref{moore-complex}.) The composition
$$
N'X_n  \stackrel{\mathrm{def}}{=} \bigcap_{i=1}^{n}\mathrm{Ker}\{d_{i}: X_{n} \longrightarrow X_{n-1} \}
\stackrel{\subseteq}{\longr} X_n \longr NX_n
$$
is an isomorphism of $R$-modules and the homomorphism
$d_0:N'X_n \to N'X_{n-1}$ yields an isomorphism of chain
complexes. We immediately drop the distinction between $NX$ and
$N'X$ and we note that Example \ref{moore-complex} implies
that $H_\ast NX$ is naturally isomorphic to the homotopy
groups of the geometric realization of $X$; hence we 
write $\pi_\ast X = H_\ast NX$.

\begin{theorem}[Dold-Kan]\label{dold-kan}The normalized chain
complex functor
$$
N: s\Mod_{R} \longrightarrow \Ch_{\ast}R
$$ 
is an equivalence of categories.
\end{theorem}

The inverse functor to $N$ is easy to write down; indeed, $N$
has a right adjoint whose definition is determined by the Yoneda Lemma.
If $C$ is chain complex, then $n$-simplices of $K(C)$ are given
by the equation
$$
K(C)_n =  \Ch_\ast R(NR[\Delta^n],C)
$$
where $R[\Delta^n]$ is the free simplicial $R$-module on the $n$--simplex. 
Using the natural isomorphism
$N'R[\Delta^n] \cong NR[\Delta^n]$ we get a formula for $K$:
\begin{equation}\label{formula-for-k}
(KC_{\bullet})_{n} = \bigoplus_{\phi:[n] \twoheadrightarrow [m]} 
\phi^\ast C_{m}
\end{equation}
where the $\phi$ are surjections in $\mathbf{\Delta}$.  For example, 
$(KC_{\bullet})_{0}=C_{0}$, $(KC_{\bullet})_{1}=C_{1} \oplus s_{0}C_{0}$, and
$$
(KC_{\bullet})_{2}=C_{2} \oplus s_{0}C_{1} \oplus s_{1}C_{1} \oplus s_{0}s_{0}C_{0}.
$$
The action of the degeneracy maps on $KC_\bullet$ is determined
by these formulas and the simplicial identities. To get the action
of the face maps $d_i$, $i > 0$, we use the simplicial identities
and require that $d_i = 0$ on $C_n \subseteq (KC_\bullet)_n$.
Finally to get the action of $d_0$, we use the simplicial identities
and require that $d_0 = \partial$ on $C_n \subseteq (KC_\bullet)_n$.

The Dold-Kan Theorem and Theorem \ref{chains-a-model} gives a model
category structure on $s\Mod_{R}$ by transport of structure. If
$X$ is a simplicial object in a category $\calC$, the underlying degeneracy
diagram neglects the face maps. More formally, if $\Delta_+$
is the subcategory of the ordinal number category $\Delta$ with
the same objects, but only surjective morphisms, then the
{\it underlying degeneracy diagram} of $X:\Delta^{op} \to \calC$
is the restriction of $X$ to $\Delta_+^{op}$.

\begin{proposition}\label{model-simp-modules} The category $s\Mod_R$
has the structure of a model category where a morphism $f:X \to Y$
is
\begin{enumerate}

\item a weak equivalence if  $\pi_{\ast}X \to \pi_{\ast}Y$ is an isomorphism;

\item  a fibration if $NX_{n} \longrightarrow NY_{n}$ is onto for $n \geq 1$; and

\item a cofibration if the underlying morphism of degeneracy diagrams
is isomorphic to a morphism of the form
$$
X_{n} \longrightarrow Y_{n} = X_{n} \oplus \bigoplus_{\phi:[n] \twoheadrightarrow [k]} \phi^\ast P_{k}
$$
where all $P_{k}$ are projective.
\end{enumerate}
\end{proposition}

It is possible to give an intrinsic description of the fibrations in
$s\Mod_R$; that is, a description that does not appeal to the
normalization functor. See Proposition \ref{fib-sim-grps}.

\begin{lemma}\label{r-mod-fibs} Let $f:X \to Y$ be a morphism
in the category $s\Mod_R$ of simplicial $R$-modules. Then the
following are equivalent.
\begin{enumerate}

\item The morphism $f$ is a fibration in $s\Mod_R$.

\item The morphism $f$ is a fibration of simplicial sets.

\item The induced map $X \to \pi_0X \times_{\pi_0Y} Y$ is a surjection.
\end{enumerate}
\end{lemma}

We next turn to the question of what it means for an $R$-module
$M$ to have a resolution in simplicial $R$-modules.
If $X$ is a simplicial object in any category, an
{\it augmentation} is a morphism $d_0:X_0 \to A$ so that
$$
d_0d_0 = d_0d_1:X_1 \longr A.
$$
In the case of simplicial modules, we automatically obtain a 
$R$-module map $\pi_0X \to A$.

\begin{definition}\label{simp-mod-res} Let $M \in \Mod_{R}$.  Then a
{\it simplicial resolution} for $M$
is an augmented simplicial $R$-module $X \to M$ such that
\begin{enumerate}
\item $X_n$ is a projective $R$-module for all $n \geq 0$; and

\item the augmentation induces an isomorphism of $R$-modules
$$
\pi_{k}X \cong \left\{ \begin{array}{cl}M & k=0 \\0 & k \neq 0. \end{array}\right.
$$
\end{enumerate}
Equivalently, a simplicial resolution for $M$ is a cofibrant replacement
for $M$ in $s\Mod_R$.
\end{definition}

Now  Theorem \ref{dold-kan} immediately yields the following
observation. Note that if $X \to M$ is an augmented simplicial
$R$-module, then $NX \to M$ is an augmented chain complex.

\begin{lemma}\label{simplicial-is-normal} Let $M$ be an
$R$-module. Then $X \to M$ is a simplicial resolution of $M$
if and only if $NX \to M$ is a projective resolution of $M$.
\end{lemma}

Now we turn to resolutions in a more general setting.
At the beginning of section 1, we developed the notion of a
resolution of an $R$-module by specifying a class of projectives
or, equivalently, a class of surjections. Then a resolution was
an acyclic chain complex of projectives; such were built by using the
fact there were enough projectives.

Now let $\Alg_{R}$ be the category of commutative $R$-algebras. 
An obvious class of projective objects in $\Alg_{R}$ are the free
commutative algebras on projective modules $P$; that is, algebras
of the form $S_{R}(P)$ where $S_R$ is the symmetric
$R$-algebra functor.\footnote{If $P = R^n$ with basis $\{x_i\}$, then
$S_R(P) \cong R[x_1,\ldots,x_n]$.} Then we should take the surjections in
$\Alg_{R}$ to be the onto maps. Resolutions will be augmented
simplicial objects in $s\Alg_R$; the following example implies
that any simplicial $R$-algebra $T$ has a natural augmentation,
in $R$-algebras, to $\pi_0Y$.

\begin{example}\label{pio-is-algebra} If $Y$ is a simplical $R$-algebra, then
$$
\pi_{0}Y = Y_{0}/(\mathrm{Im}(d_{0} - d_{1})).
$$
But, because of the presence of the degeneracy morphism $s_0$, the
$R$-module  Im$(d_{0}-d_{1})$ is an ideal in $Y_0$: if $a \in Y_{0}$
and $x=(d_{0}-d_{1})(y)$, then
\begin{align*}
(d_{0}-d_{1})(s_{0}(a)y)  &= d_{0}s_{0}(a)\cdot d_{0}(y) - d_{1}s_{0}(a)d_{1}(y) \\
& =  a(d_{0}(y) - d_{1}(y))\\
&= a[(d_{0}-d_{1})(y)]
\end{align*} 
since the $d_i$ are algebra morphisms and $d_{i+1}s_{i} = d_{i}s_{i} = 1$.
\end{example}

Next, a bit of notation.
Let $X$ be a simplicial object in any category. We will write $X_+$ for the
underlying degeneracy diagram; that is, if $\ddelta_+$
is the subcategory of the ordinal number category where the morphisms
are the surjections, then $X_+$ is $X$ restricted to $\ddelta^{op}_+$.

\begin{definition}\label{free-object-resolution} A {\it resolution} of
$A \in \Alg_{R}$ is an augmented simplicial $R$-algebra $X \to A$
so that:
\begin{enumerate}
\item there is a functor $P:\ddelta^{op}_+ \to \Mod_R$ so that
$P_n$ is projective for all  $n$ and so that there is an isomorphism
of diagrams of $R$-algebras $X_+ \cong S_R(P)$

\def\defeq{{\stackrel{def}{=}}}

\item the augmentation induces an isomorphism
\[
\pi_{k}X \cong \left\{ \begin{array}{cl} A & k = 0 \\ 0 & k \neq 0 \end{array} \right.
\]
\end{enumerate}
\end{definition}

\begin{remark}\label{comments-on-simp-res} 1.) In point (1) of
this definition,we are requiring more than that $X$ be level-wise projective; we are requiring that the degeneracy maps respect this
fact is a systematic way. This should be compared to Definition
\ref{simp-mod-res}.1 where we do not seem to make any such
requirement. However, because $R$-modules form an abelian
category, the analogous requirement would be automatic. 

2.) We have said nothing about
the existence or
uniqueness of such resolutions, but we will see that such a resolution
turns out to be a cofibrant replacement for $A$ in an appropriate
model category structure on simplicial algebras, so Corollary
\ref{uniqueness-of-replacements} applies. See Proposition
\ref{free-cof}.
\end{remark}

\begin{example}\label{extra-degeneracy} We will see that the
adjoint pair
$$
\xymatrix{
S_R:s\Mod_R \ar@<.5ex>[r]& \ar@<.5ex>[l]\ s\Alg_R:\cO
}
$$
(where $\cO$ is the forgetful functor) will be a Quillen functor. This
implies that $S_R$ preserves weak equivalences between cofibrant
objects, but in this example we will be much more concrete in hopes
of further explaining why it is worthwhile to use simplicial objects.
These observations are directly from Dold's paper \cite{dold}.

In any category, an augmented simplicial object $X \to X_{-1}$ has an
{\it extra degeneracy} if there are morphisms $s_{-1}:X_n \to X_{n+1}$
$n \geq -1$ so that the
formulas of Lemma \ref{simp-relations} still hold. This implies
that if $X$ is a simplicial $R$-module, then
the normalized augmented chain complex $NX$ has a chain contraction
onto $X_{-1}$ given by setting $T = s_{-1}:NX_n \to NX_{n+1}$.
The Dold-Kan
correspondence implies that $X$ has an extra degeneracy if
and only if $NX$ has a chain contraction. 

Now, if $X$ is a simplicial $R$-module with an extra degeneracy,
then $S_R(X)$ has an extra degeneracy and, hence
$\pi_\ast S_R(X) \cong S_R(X_{-1})$. In particular, if $X$ is
cofibrant and $\pi_\ast X = 0$, then $\pi_\ast S_R(X) \cong
R$ in degree $0$. By contrast, if $C_\bullet$ is a chain
complex with a chain contraction, the associated differential
graded algebra $S_R(C_\bullet)$ does not inherit a contraction.
Again, see Example \ref{something-to-prove}.

The notion of an extra degeneracy can be generalized to the notion
of a simplicial homotopy (see \cite{May}) and we are using the
observation that the functor $S_R$ preserves simplicial homotopies.
\end{example}

\subsection{Simplicial model categories}

In the previous section we tried to suggest that simplicial resolutions
were an appropriate and natural generalization of projective resolutions
for $R$-modules. In turns out, in addition, that categories of simplicial
objects support a great deal of structure -- in particular, they are
enriched over simplicial sets -- and it is extremely useful to have
an interplay between this structure and any model category structure 
that may be around. This leads to the notion of a simplicial
model category.

If $\calC$ is a category, let $s\calC$ denote the simplicial objects in
$\calC$. We will assume that $\calC$ has all limits and colimits.

We first note that $s\calC$ has an {\it action} by simplicial sets. If $K \in
\sisets$ and $X \in s\mathcal{C}$, define $K \otimes X \in s\calC$ as
the coproduct
\begin{equation}\label{k-tensor-x}
(K \otimes X)_{n} = \coprod_{K_{n}} X_{n}.
\end{equation}
The face and degeneracy maps are determined by those in $K$ and
$X$. This construction gives  a bifunctor $\sisets \times s\mathcal{C} \longrightarrow s\mathcal{C}$ and one easily checks that there are
natural isomorphisms 
\begin{equation}\label{sm01}
(K \times L) \otimes X \cong K \otimes(L \otimes X)
\end{equation}
and
\begin{equation}\label{sm02}
{\Delta}^{0} \otimes X = \ast \otimes X \cong X.
\end{equation}
Less obvious, but equally formal is that there is also an {\it exponential}
bifunctor
\begin{align*}
s\calC \times \sisets^{op} &\longr \calC\\
(Y,K) &\longmapsto Y^K
\end{align*}
determined by the adjoint formula
$$
s\calC(K \otimes X,Y) \cong s\calC(X,Y^K).
$$

\def\map{{{\mathrm{map}}}}

From this data we construct a simplicial mapping space.  Let $X,Y \in
s\mathcal{C}$, define $\mathrm{map}_{s\mathcal{C}}(X,Y)_{n} \in
\sisets$ by
\[
\mathrm{map}_{s\mathcal{C}}(X,Y)_{n} = s\calC (\Delta^n \otimes
X, Y).
\]
A morphism $\phi$ in the ordinal number category determines a
morphism $\phi_\ast:\Delta^m \to \Delta^n$ which, in turn, gives
$\map_{s\calC}(X,Y)$ its structure as a simplicial set. The formulas
of Equations \ref{sm01} and \ref{sm02} supply 
 an associative composition
\[
\mathrm{map}_{s\mathcal{C}}(Y,Z) \times \mathrm{map}_{s\mathcal{C}}(X,Y) \longrightarrow \mathrm{map}_{s\mathcal{C}}(X,Z).
\]
and a natural isomorphism
$$
\map_{s\calC}(X,Y)_{0}=s\calC(X,Y).
$$
Thus we have shown that $s\mathcal{C}$ is {\it enriched}
over simplicial sets. But more is true: there are enriched
adjoint isomorphisms
$$
\map_{s\calC}(K \otimes X,Y) \cong \map_{s\calC}(X,Y^K) \cong
\map_{\sisets}(K,\map_{s\calC}(X,Y)).
$$
All this data together gives $s\calC$ the structure of a {\it simplicial
category}, as defined by Quillen in \cite{HA} \S II.2.

\begin{examples}[Examples of simplicial categories]\label{exams-simp}
1. If $s\calC = \sisets$, then $K \otimes X =
K \times X$ and
$$
X^K = \map_{\sisets}(K,X).
$$
This example predates Quillen, of course. 

2. If $s\calC = s\Mod_R$, then
$$
K \otimes X = R[K] \otimes_R X
$$
where $R[-]$ is the free $R$-module functor extended level-wise to 
simplicial sets. In addition
$$
X^K = \map_{\sisets}(K,X)
$$
with addition and $R$-module action coming from the target. Similarly,
if $X \in s\Alg_R$ where $\Alg_R$ is the category of commutative $R$-algebras, then
$X^K = \map_{\sisets}(K,X)$ with algebra structure arising from the
target. In this case $K \otimes X$ must be calculated using Equation
\ref{k-tensor-x}. But this is not so bad, the coproduct in commutative
$R$-algebras is given by tensor product.

3. This last example can be greatly expanded to almost any reasonable
algebraic structure. Thus, to name just a few cases, there is a simplicial
model category structure on simplicial associative algebras, simplicial
Lie algebras, simplicial restricted Lie algebras over a field of
positive characteristic, simplicial algebras over any operad in sets,
and simplicial unstable algebras over the Steenrod algebra. All
of these examples have appeared in the literature.

4. The category $\cgh$ of compactly generated weak Hausdorff
spaces is a simplicial category. If $K$ is a simplicial set and
$X \in \cgh$, set $K \otimes X = |K| \times X$ and $X^K = X^{|K|}$,
where the $X^{|K|}$ denotes the exponential object in $\cgh$.
The reader sensitive to category theory will point out that
$\cgh$ is really a {\it topological} category -- defined by analogy
with simplicial category -- to which we respond that every topological
category is a simplicial category via geometric realization.\footnote{This
is a consequence of the fact that geometric realization preserves
finite products. While this point is elementary in the sense that it is pure
category theory, it does have its subtleties.}

5. Any simplicial category has a ready-made and natural
candidate for a path object; namely, $X^{\Delta^1}$. Path objects
arose in Theorem \ref{quillen-criterion}.
\end{examples}

We would now like to meld the notion of simplicial structure with the 
structure of a model category when it is present. This is encoded
in the following axiom, which in the literature is known as Quillen's
SM7.

\begin{definition}[Corner Axiom]\label{corner-axiom}
Let $\mathcal{C}$ be a category which is at once a model category and
a  simplicial category. Then $\mathcal{C}$ is a {\it simplicial model 
category} if for $j: A \longrightarrow B$ a cofibration and $q: X \longrightarrow Y$ a fibration, the natural map of simplicial mapping
spaces
\[
\mathrm{map}_{\calC}(B,X) \longr
\mathrm{map}_{\calC}(B,Y) \times_{\mathrm{map}_{\calC}(A,Y)}
\mathrm{map}_{\calC}(A,X)
\]
is a fibration of simplicial sets which is a weak equivalence if
$j$ or $q$ is a weak equivalence in $\calC$.
\end{definition}

The corner axiom can be reformulated in terms of the action of
simplicial sets on $\calC$ -- often giving a condition which is easier
to check. The proof is an easy exercise in adjointness arguments.

\begin{proposition}\label{other-sm7s}Let $\calC$ be a model
category. The corner axiom is equivalent
to both of the following two statements. 

1. Let $j: A \to B$ be a cofibration in $\calC$ and
$i:K \to L$ be a cofibration of simplicial sets.  Then
\[
i \otimes j: L \otimes A \sqcup_{K \otimes A} K \otimes B \longrightarrow
 L\otimes B
\]
is a cofibration in $\mathcal{C}$ which is a weak equivalence if $i$ or $j$ is.

2. Let $q:X \to Y$ be a fibration in $\calC$ and $i:K \to L$ a cofibration
of simplicial sets. Then
$$
X^L \longr X^K \times_{Y^K} Y^L
$$
is a fibration in $\calC$ which is a weak equivalence if $q$ or $i$ is.
\end{proposition}

If $\calC$ happens to be cofibrantly generated, there can be a further
reduction in Proposition \ref{other-sm7s}.1: we need only check
the condition on $i \otimes j$ where $i$ and $j$ run over the
generating  cofibrations (or acyclic cofibrations, as needed) for
$\calC$ and for $\sisets$. This is because, being a left adjoint in
both variables $(-) \otimes (-)$ commutes with all colimits. These
observations make it relatively easy to prove that there
are examples of simplicial model categories:

\begin{theorem}\label{exams-simp-model} The categories $\sisets$,
$s\Mod_R$, and $\cgh$, with the simplicial structures of
Example \ref{exams-simp} become simplicial model categories.
\end{theorem}

We close this section with two remarks meant to underscore the
fact that simplicial model categories have the right sort of rich
structure for us.

\begin{corollary}\label{path-cylinder}Let $\calC$ be a simplicial
model category.
\begin{enumerate}

\item If $X \in \mathcal{C}$ is cofibrant, then 
$$
X \amalg X \cong \Delta^{1} \otimes X  \longr
\partial \Delta^{1} \otimes X \longr \ast \otimes X \cong X
$$
is a natural cylinder object for $X$.

\item If $X$ in $\calC$ is fibrant, then 
$$
X \longr X^{\Delta^1} \longr X^{\partial\Delta^1} \cong X \times X
$$
is a natural path object for $X$.
\end{enumerate}
\end{corollary}

The next result says that the mapping space is a rich homotopical
object. Indeed, one of the lessons of the last thirty years or so is
that in order to compute homotopy classes of maps, the best strategy
can be to compute the homotopy type of the mapping space, then
read off the components.

\begin{corollary}\label{mapping-homotopy}Let $\calC$ be a simplicial
model category.
If $X \in \mathcal{C}$ is cofibrant and $Y \in s\mathcal{C}$ is fibrant, then
\[
\pi_{0}\mathrm{map}_{\mathcal{C}}(X,Y) = \Ho(\mathcal{C})(X,Y)
= [X,Y]_{\mathcal{C}}.
\]
\end{corollary}

\begin{remark}\label{hammock} In a series of important papers, Dwyer and
Kan noted that given {\it any} model category $\calC$, there is
 a mapping simplicial set $L_H\calC(X,Y)$ between two objects $X$ and $Y$ of
$\calC$ with the properties that 
\begin{enumerate}
\item $\pi_0L_H(X,Y) = [X,Y]_\calC$, and

\item if $\calC$ is a simplicial model category, $X$ is cofibrant,
and $Y$ is fibrant,
there is a natural zig-zag of weak equivalences from $\map_\calC(X,Y)$
to $L_H(X,Y)$.
\end{enumerate}

This is the Dwyer-Kan {\it hammock localization} of $\calC$. They point out, in
fact, that one really only needs the weak equivalences to define $L_H(X,Y)$.
See \cite{hammock}, \cite{classification}, and \cite{DKTop}.
Of course, the simplicial mapping space $\map_\calC(X,Y)$ is
more concrete and the corner axiom of Definition \ref{corner-axiom} and 
its adjunct Proposition \ref{other-sm7s} are very useful for computations.
\end{remark}

\subsection{Simplicial algebras}

We'd now like to pull together the threads from cofibrantly
generated model categories and simplicial model categories.
Here is an example of the sort of result we'd like to prove.
The model category on simplicial $R$-modules was discussed
in Proposition \ref{model-simp-modules}. The following result
leaves open a characterization of cofibrations. While they will
be formally determined -- see Remark \ref{model-cat-remarks}.2 --
we want to be  more concrete. This will be addressed below.

\begin{theorem}\label{simp-R-alg-model}
Let $R$ be a commutative ring and $s\Alg_R$ the
category of simplicial commutative algebras over $R$. Then 
$s\Alg_R$ has the structure of a simplicial model category
where a morphism $f:X \to Y$ is
\begin{enumerate}

\item a weak equivalence if $\pi_\ast X \to \pi_\ast Y$ is an
isomorphism and

\item a fibration if the induced map $X \to \pi_0X \times_{\pi_0Y} Y$
is a surjection.
\end{enumerate}
\end{theorem}

\begin{proof}We would like to apply Theorem \ref{lifting} to the
adjoint pair
$$
\xymatrix{
S_R:s\Mod_{k} \ar@<2pt> [r] & s\Alg_{k}:\cO \ar@<2pt> [l]
}
$$
where $\cO$ is the forgetful functor and $S_R$ is the symmetric
algebra functor. Note that, by Lemma \ref{r-mod-fibs}, we have
{\it defined} a morphism $f\in s\Alg_R$ to be a  weak equivalence or
fibration if and only if $\cO(f)$ is a weak equivalence or fibration.
Since $S(-)$ commutes with all filtered colimits, we can appeal
to Theorem \ref{quillen-criterion} to complete the model category
structure. To check that result, we turn to Corollary \ref{path-cylinder}.
If $B \in s\Alg_R$,
then $B^{\Delta^1}$ is a path object for $B$ in $s\Mod_R$ --
because $s\Mod_R$ is already a simplicial model category. It
follows immediately that $B^{\Delta^1}$ is a path object in
$s\Alg_R$. Since {\it every} object of $s\Alg_R$ is fibrant,
Theorem \ref{quillen-criterion} applies, and we have our model
category structure.

To a get a {\it simplicial} model category structure, we appeal
to Proposition \ref{other-sm7s}.2. Again we note that it is sufficient
to check the condition there in $s\Mod_R$.
\end{proof}

\begin{remark}\label{so-formal}The reader will have remarked that the
previous argument is very formal and has wide application. In particular,
we have simplicial model category structures on simplicial associative
algebras, simplicial Lie algebras, simplicial groups, and so on.
\end{remark}

The next point is to characterize cofibrations in $s\Alg_R$. Again,
Theorem \ref{lifting} can help, and  this gives a chance to introduce
a new construction and some other new ideas.

Let $I$ be a small category, $\calC$ any category with colimits
and $\calC^I$ the category of $I$-diagrams in $\calC$. Let $I^\delta$ be
the category with same objects as $I$ but only identity morphisms; thus
$I^\delta$ is $I$ made discrete. An $I$-diagram $X:I \to \calC$ is
$I$-{\it free} (or simply free) if it is the left Kan extension of
some diagram $Z:I^\delta \to \calC$. In formulas, this means
that 
$$
X_i = \mathop{\sqcup}_{j \to i} Z_j
$$
with the coproduct over all morphisms $j \to i$ in $I$.

\begin{definition}\label{s-free} Let $\ddelta$ be the ordinal number category and
$\ddelta_+ \subset \ddelta$ the category with same objects but only surjective
morphisms. Let $\calC$ be a category and $X:\ddelta^{op} \to \calC$
a simplicial object. Then $X$ is $s$-{\it free} if the underlying diagram
$$
X_+:\ddelta_+^{op} \longr \calC
$$
is free. More generally, a morphism $X \to Y$ is of simplicial objects
is $s$-free if the underlying morphism $X_+ \to Y_+$ of $\ddelta_+^{op}$
diagrams is isomorphic to the inclusion of a summand $X_+ \to X_+ \amalg Y_0$
where $Y_0$ is $s$-free.
\end{definition} 

To rephrase, to say that $X$ is $s$-free is to say that there are objects
$Z_k$ so that there are isomorphisms
$$
X_n = \coprod_{\phi:[n] \to [k]} \phi^\ast Z_k
$$
where $\phi$ runs over the surjections in $\Delta$. In this case we
say $X$ is $s$-free on $\{Z_k\}$. More generally, a morphism can
be $s$-free on a set of objects.

\begin{definition} A morphism $X \to Y$ of simplicial $R$-algebras
is {\it free}\footnote{The terminology, while not fabulous, is
in \cite{HA} and, thus, hallowed by history.} if it is $s$-free on
a set of objects $\{S_R(P_k)\}$ where each $P_k$ is a projective
$R$-module.
\end{definition}

The following result characterizes cofibrations in $s\Alg_R$. It also
shows that if we regard an $R$-algebra $A$ as a constant simplicial
$R$-algebra, then it has a resolution in the sense of 
Definition \ref{free-object-resolution}.

\begin{proposition}\label{free-cof}A morphism $X \to Y$ in $s\Alg_R$ is a cofibration
if and only if it is a retract of a free morphism. Furthermore, any
morphism can be factored naturally as 
$$
\xymatrix{
X \rto^i & Z \rto^p & Y
}
$$
where $i$ is a free map and $p$ is an acyclic fibration. A similar
statement applies to the ``acyclic cofibration-fibration'' factorization
as well.
\end{proposition}

\begin{proof} Note that free maps are closed under coproducts, cobase
change, and sequential colimits. Since the generating cofibrations and
the generating acyclic cofibrations are free maps, the statement about
the factorizations follows from the small objects argument. Then
if $f:A \to B$ is any cofibration, we obtain a diagram and a lifting
problem
$$
\xymatrix{
A \rto^i \dto_f & Z \dto^p \\
B \rto_{=} \ar@{-->}[ur] & B
}
$$
where $i$ is a free map and $p$ is an acyclic fibration. Since $f$ is
a cofibration, the lifting problem can be solved and we have that
$f$ is a retract of $i$. Hence every cofibration is a retract of a free
map. It remains to be shown that every free map is a cofibration.
This will be done below, when we have the skeletal decomposition.
See Lemma \ref{free-is-cob}.
\end{proof}

\subsection{Homology and Cohomology}

To illustrate the efficacy of model categories, we will give
Quillen's definition of the homology of commutative algebras.
The paper \cite{AQ} is still a wonderful read.
This was one of the first real applications of model categories,
and the proof of the transitivity sequence (Proposition \ref{trans})
is most easily given using all the language of the theory.

\begin{definition}\label{abelian-objects}
If $\mathcal{C}$ is a category, $A \in \mathcal{C}$ is an {\it abelian object}
if $\calC(-, A)$ is naturally an abelian group.  
Assuming $\calC$ has enough limits, this is equivalent to there being a multiplication morphism
$$
m: A \times A \longrightarrow A
$$
an identity $\epsilon: \ast \longrightarrow A$, where $\ast$ is the terminal object, and an inverse $i: A \longrightarrow A$ such that all the usual diagrams commute.
\end{definition}

We will let $\mathcal{C}_{\mathit{ab}}$ denote the subcategory of abelian objects in $\mathcal{C}$.

\begin{examples}\label{exams-abelian}
1.) If $\mathcal{C}$ is the category of sets, then
$\mathcal{C}_{\mathit{ab}}$ is
the category of abelian groups.  If $\mathcal{C} = \sisets$, then
$\mathcal{C}_{\mathit{ab}} = s\Mod_{\mathbb{Z}}$.
\smallskip

2.) Fix a commutative ring $R$. The only abelian object in the category 
$\Alg_R$ of commutative $R$-algebras is the terminal obect $0$. To
get more interesting examples we work over a fixed $R$-algebra $A$.  
Then if $M$ is an $A$-module, define a new $R$-algebra over
$A$ by setting  $A \ltimes M = A \oplus M$ with the multiplication
$$
(a,x)(b,y) = (ab, ay + bx).
$$
 Then
\begin{eqnarray*}
(\Alg_{R}/A)(X, A \ltimes M) & \cong & \mathrm{Der}_{R}(X, M) \\
f = \epsilon \oplus \partial & \longmapsto & \partial.
\end{eqnarray*}
Here $\mathrm{Der}_{R}(X, M)$ is the $R$-module of derivations;
that is, of $R$-module maps $X \to M$ satisfying the Leibniz rule.
This equation displays $A  \ltimes M$ as an abelian object.  In fact,
a simple exercise, shows that  every abelian object in $\Alg_R/A$ is of
this form.  In particular, the functor $A\ltimes (-)$ defines
equivalences of categories
$$\Mod_A \to (\Alg_{k}/A)_{ab}$$
and
$$
s\Mod_{A} \longr (s\Alg_{k}/A)_{ab}.
$$
\end{examples}

For the next definition, we will suppose there are model category structures
on $\calC$ and $\calC_\mathit{ab}$ so that the inclusion
$\mathcal{C}_{ab} \longrightarrow \mathcal{C}$ is the right adjoint of
a Quillen functor; that is, there is an abelianization functor 
$$
\Ab: \mathcal{C} \longrightarrow \mathcal{C}_{ab}
$$
left adjoint to inclusion which preserves cofibrations and weak equivalences
between cofibrant objects.

\begin{definition}\label{homology}Homology is the total left derived
functor of abelianization. That is, if $X \in \calC$, then the {\it Quillen
homology} of $X$ is the object $L\Ab(X) \in \calC_\mathit{ab}$.
\end{definition}

\begin{example}[Homology of spaces]\label{sing-homology}
If $X \in \sisets$, then $X$ is cofibrant and $L\Ab(X) = \mathbb{Z}X$, the simplicial abelian group generated by $X$.
In particular,  we have that if $Y$ is a topological space
$$
\pi_{n}L\Ab(S(Y)) = \pi_n\mathbb{Z}S(Y) \cong H_{n}(Y; \mathbb{Z})
$$
and we recover singular homology.

If $X$ is a topological space, then a cofibrant replacement for $X$ can be taken
to be a weak equivalence $Y \to X$ with $Y$ cofibrant. Then 
$\Ab(Y)$ is the free topological abelian group on $Y$ and the Dold-Thom
Theorem implies that there are natural isomorphisms
$$
\pi_\ast L\Ab(X) \cong \pi_\ast\Ab(Y) \cong H_\ast(Y) \cong H_\ast (X).
$$
\end{example}

\begin{example}[Homology of groups]\label{group-homology}If
$G$ is a group, we may regard $G$ as a constant simplicial group,
and then $L\Ab(G) = \Ab(X) = X/[X,X]$ where $X \to G$ is
a cofibrant model for $G$ and $[X,X]$ is the commutator subgroup.
On the other hand the homology of $G$ is usually defined to be
$$
\mathrm{Tor}^{\ZZ[G]}_\ast (\ZZ,\ZZ) = H_\ast (BG)
$$
where $BG$ is the classifying space of $G$. The claim is that there is
a degree-shifting isomorphism between $\pi_\ast L\Ab(G)$ and
$\tilde{H}_\ast (BG)$. To see this, form the bisimplicial abelian group 
$\ZZ[BX]$. This has two filtrations, and hence, two spectral sequences
converging to the same graded abelian group.

For all $q$, $X_q$ is a free group. Thus,
if we filter by the simplicial degree coming from $X$ we get a
spectral sequence with $E^1$-term 
$$
E^1_{p,q} \cong H_q(BX_p) \cong \left\{ \begin{array}{ll}
                                {X_p/[X_p,X_p]}&{q=1}\\
                                {\ZZ}&{q=0}\\
				    {0}&{q > 1}\end{array}\right.
$$
Thus $E^2_{p,1} = \pi_pL\Ab(G)$, $E^2_{0,0} = \ZZ$ and
$E^2_{p,q} = 0$ for all other $p$ and $q$. One the other
hand, if we filter by the simplicial degree coming from the
functor $B(-)$ we get a spectral sequence with
$$
E^1_{p,q} \cong \left\{ \begin{array}{ll}
                                {\ZZ[BG]}&{p=0}\\
                                \\
                                {0}&{p \ne 0}\end{array}\right.
$$
Thus $E^2_{0,q} = H_qBG$ and $E^2_{p,q} = 0$ if $p \ne 0$.
Thus we conclude that
$$
\pi_nL\Ab(G) \cong H_{n+1}BG.
$$
\end{example}

We now turn to the case of commutative algebras over a commutative
ring $R$. The resulting homology theory is {\it Andr\'e-Quillen
homology} and it is the subject of a very extensive discussion
elsewhere in these nots by Srikanth Iyengar \cite{Iyengar}, which the reader should
turn to. We only include it here because
(a) the construction of a suitable homology theory for commutative
algerbas was one of the early successes of the theory of model categories
and (b) the first author of this monograph is very fond of it.

The first question is what the abelianization functor for commutative
algebras should be. 

If $X$ is an $R$-algebra, let $I = \mathrm{Ker}(X \otimes_R X \to
X)$ be the kernel of $R$-algebra multiplication and let
$$
\Omega_{X/R} = I/I^2.
$$
Then the morphism $d:X \to \Omega_{X/R}$ sending $y$ to
the coset of $y \otimes 1 - 1\otimes y$ is a derivation
and the natural map
\begin{align*}
\Mod_X(\Omega_{X/R},M) &\longr \Der_R(X,M)\\
f&\mapsto f \circ d
\end{align*}
is an isomorphism. Since $\Omega_{X/R}$ represents derivations, it is
called the module of differentials.

Let $A$ be an $R$-algebra and $\Alg_R/A$ the category of algebras
over $A$.  Recall that an object in $\Alg_R/A$ is a morphism of 
commutative $R$-algebras $B \to A$. In Example \ref{exams-abelian}.2 we
indicated that the
category of abelian objects in $\Alg_R/A$ was equivalent to the 
category $\Mod_A$ and that the inclusion functor of abelian
objects into $\Alg_R/A$ was naturally isomorphic to the functor
$$
A \ltimes (-): \Mod_R \longr \Alg_R/A.
$$
But there are natural isomorphisms
$$
(\Alg_R/A)(X,A\ltimes M) \cong \Der_R(X,M) \cong \Mod_A(A 
\otimes_X \Omega_{X/R},M).
$$
Therefore, we have proved the following result, which displays
the differentials as the abelianization functor for $\Alg_R/A$.

\begin{proposition}\label{differentials-abelian}The functor
$$
X \mapsto A \otimes_X \Omega_{X/R}
$$
is left adjoint to the functor $M \mapsto A \ltimes M$ from
$\Mod_A$ to $\Alg_R/A$.
\end{proposition}

Note that the functor $M \mapsto A \ltimes M$ extends level-wise
to a functor $s\Mod_A \to s\Alg_R/A$. This functor preserves
all weak equivalences and fibrations, so the adjoint
$$
\Omega_{(-)/R}:s\Alg_R/A \longr s\Mod_A
$$
preserves cofibrations and weak equivalences between cofibrant objects;
that is, we have a Quillen functor.

\begin{definition}\label{aq-homology}Let $A$ be a commutative
$R$-algebra. Then the {\it cotangent
complex} for $A$ is the simplicial $A$-module given by the total left
derived functor of differentials:
\[
L_{A/R} =  A \otimes_{X} \Omega_{X/R},
\]
where $X \to A$ is a cofibrant replacement for $A$ in
the category $s\Alg_R/A$. (Cofibrant objects in $s\Alg_R/A$ can be understood
using Proposition \ref{free-cof} and Example \ref{formal}.) The
{\it Andr\'e-Quillen} homology of $A$ is given by
$$
D_q(A/R) = \pi_qL_{A/R} = H_qNL_{A/R}.
$$  
More generally, if $M$ is an $A$-module, we set $D_\ast(A/R;M) = \pi_\ast 
(M \otimes_A L_{A/R})$.
\end{definition}

\begin{remark}\label{nat-of-aq}We write down here that 
there is a natural augmentation $L_{A/R} \to \Omega_{A/R}$.

It may not be immediately 
obvious what the naturality properties of $L_{A/R}$ should be.
The important observation is that the cotangent complex
is a {\it relative} homology object; that is, a functor of
the arrow $R \to A$ -- or, we might say, of the pair $(A,R)$.
If we have a commutative diagram
$$
\xymatrix{
R \dto\rto & S\dto\\
A \rto & B
}
$$
then we get a morphism $B\otimes_A L_{A/R} \to L_{B/S}$.
To see this, form the following diagram, which can be chosen
naturally: 
$$
\xymatrix{
X \rto \dto_p & S \otimes_R X \dto \rto^-i & Y \dto^q\\
A \rto & B \rto^{=} &B.\\
}
$$
Here $p$ is an acyclic fibration in $R$-algebras with $X$
cofibrant, $i$ is a cofibration in $S$ algebras and $q$ 
is an acyclic fibration in $S$-algebras. Then $S \otimes_RX$
is a cofibrant $S$-algebra and, hence, $Y$ is a cofibrant
$S$-algebra. Thus we may model $B\otimes_A L_{A/R} \to L_{B/S}$
by
$$
B \otimes_A (A \otimes_X \Omega_{X/R}) \to B \otimes_Y
\Omega_{Y/S}.
$$
\end{remark}

We might ask to what extent this sort of homology actually
acts like homology. We'll give two answers. First, the next
result gives some basic properties of this theory, including
a result that says Andr\'e-Quillen homology takes finite coproducts
to sums. (The finite coproduct of commutative algebras is their
tensor product.) Second, we'll also give below, in Theorem \ref{trans}
a proof of the transitivity sequence, which is an analog of the
long exact sequence of a pair.

\begin{lemma}\label{basic-props-aq}The cotangent complex has the following properties:
\begin{enumerate}

\item if $A = S_R(P)$ where $P$ is a projective $R$-module than the
natural augmentation $L_{A/R} \to \Omega_{A/R}$ is a weak equivalence
of simplicial $A$-modules;

\item if $A$ and $B$ are two $R$-algebras, let $C = A \otimes_R B$.
If one of $A$ or $B$ is flat over $R$, then the natural map
$$
C \otimes_A L_{A/R} \oplus C \otimes_B L_{B/R} \longr
L_{C/R}
$$
is an isomorphism.
\end{enumerate}
\end{lemma}

\begin{proof} The proofs are easy and we include them only to 
illustrate the flexibility of model categories.
The first statement follows from the fact that
$A = S_R(P)$ is already cofibrant in $\Alg_R/A$. For the
second statement, let $X \to A$ and $Y \to B$ be
cofibrant replacements for $A$ and $B$ respectively.
Then $X \otimes_R Y \to A \otimes_R B= C$ is a cofibrant
replacement for $C$. To see this, note that (1) the coproduct of
cofibrant objects is cofibrant and that (2) since $X$ is level-wise projective
as an $R$-module, that flatness hypothesis implies
$$
\pi_\ast (X \otimes_R Y) = \pi_\ast X \otimes_R \pi_\ast Y \cong
A \otimes_R B.
$$
The result now follows from the formula, easily checked, that
$$
\Omega_{(X \otimes_R Y)/R} = (X \otimes_R Y) \otimes_X
\Omega_{X/R} \oplus (X \otimes_R Y) \otimes_Y \Omega_{Y/R}.
$$
\end{proof}

A fundamental result about Andr\'e-Quillen homology is the
following, which is also easy to prove using model categories:

\begin{theorem}[Flat base change]\label{flat-base-change}Suppose
there is a commutative diagram
$$
\xymatrix{
R \dto_i\rto^f & S\dto\\
A \rto & S \otimes_R A
}
$$
with one of $i$ or $f$ flat. Let $B = S\otimes_R A$. Then the induced morphism of simplicial $B$-modules $B \otimes_A L_{A/R} \to L_{B/S}$
is a weak equivalence.
\end{theorem}

\begin{proof} Choose a cofibrant model $X \to A$ for $A$ as
an $R$-algebra. Then $S \otimes_R X \to B$ is a weak equivalence,
by the flatness assumption, and $S \otimes_R X$ is cofibrant.
The result follows.
\end{proof}

To set the stage for the next result, we consider a sequence of morphisms
$$
\xymatrix{
R \rto & A \rto^f & B
}
$$
of commutative rings. Then a standard result (and an easy
exercise) about differentials is that there is an exact sequence
\begin{equation}\label{proto-trans}
B \otimes_A\Omega_{A/R} \longr \Omega_{B/R} \longr \Omega_{B/A} \to 0
\end{equation}
and that, furthermore, this sequence becomes short exact if $f:A \to B$
has an $R$-algebra retraction; that is, if there is an $R$-algebra map $r:B \to A$
so that $rf = 1:A \to A$. In trying to develop a homology theory of commutative
algebras, the question was how to extend this exact sequence to the left. 
To state a result, call a sequence of simplicial $R$-modules $M \to N \to P$ 
a {\it cofiber sequence}, if the composite is null-homotopic and
if the mapping cone of the first map is weakly equivalent to $P$. Such a sequence
induces a long exact sequence in homotopy.

\begin{proposition}\label{trans}Let $R \to A \to B$ be a sequence of commutative
rings. Then the composition
$$
B \otimes_A L_{A/R} \longr L_{B/R} \longr L_{B/A}
$$
is a cofiber sequence of simplicial $B$-modules. In particular, if $M$ is a
$B$-module there is a long exact
sequence
\begin{align*}
\cdots \to D_1(A/R;M) \to D_1(B/R;M) &\to D_1(B/A;M)\\
\to M \otimes_A \Omega_{A/R}
&\to M \otimes_B \Omega_{B/R} \to M \otimes_B \Omega_{B/A} \to 0.
\end{align*}
\end{proposition}

\begin{proof}The long exact sequence follows from the cofiber sequence, as
$L_{B/A}$ will be level-wise a projective $B$-module. To prove that we have
a cofiber sequence, choose an acyclic fibration $X \to A$ where
$X$ is free. Then factor the composition $X \to A \to B$ to obtain a commutative
square
$$
\xymatrix{
X \rto^i \dto & Y \dto^p\\
A \rto& B
}
$$
where $i$ is a free map and $p$ is an acyclic fibration. Then we have an exact
sequence of simplicial $B$-modules
$$
0 \to B \otimes_X \Omega_{X/R} \to B \otimes_Y \Omega_{Y/R} \to B \otimes_Y
\Omega_{Y/X} \to 0
$$
using that Equation \ref{proto-trans} becomes short exact because
$X_n \to Y_n$ has a retract for all $n$; indeed $X_n \to Y_n$ is isomorphic to
\begin{equation}\label{display}
X_n \longr X_n \otimes_R S_R(P_n)
\end{equation}
where $P_n$ is a projective $R$-module. To finish the argument, note that
we have a map $A \otimes_X Y \to B$. This is a weak equivalence and
$A \otimes_X Y$ is a cofibrant simplicial $A$-algebra. Thus we need only 
note that Equation \ref{display} implies that 
$$
B \otimes_Y \Omega_{Y/X} \to B \otimes_{(A \otimes_X Y)}
\Omega_{(A \otimes_X Y)/A}
$$
is an isomorphism.
\end{proof}

\begin{remark}\label{etale}Andr\'e-Quillen homology is a very
sensitive invariant of morphims of commutative rings or, more
generally, of schemes. (See \cite{Ill} for the scheme-theoretic
generalizations.) For example, suppose $A$ is a finitely generated
$R$-algebra. Then $R \to A$ is smooth if and only if
$L_{A/R} \to \Omega_{A/R}$ is a weak equivalence and $\Omega_{A/R}$
is projective as an $A$-module. In addition, $R\to A$ is \'etale if and
only if $L_{A/R}$ has vanishing homology.
\end{remark}

\begin{remark}\label{cohomology}There is a companion Andr\'e-Quillen
cohomology theory. If $M$ is an $A$-module, let $K(M,n) \in s\Mod_A$
be the unique simplicial $A$-module so that
\[
NK(M,t)_{n} = \left\{ \begin{array}{cl} M, & n = t; \\ 0, & n \neq t. \end{array} \right.
\]
Then
\[
\pi_{n}K(M,t) \cong H_{n}NK(M,t) = \left\{ \begin{array}{cl} M, & n = t; \\ 0, & n \neq t.\end{array} \right.
\]
The $n$-th Andr\'e-Quilen cohomology group of $A$ is defined to be homotopy
classes of maps into the resulting Eilenberg-MacLane algebra:
\begin{align*}
D^n(A/R;M) &= \Ho(s\Alg_{R}/A)(A, A \ltimes K(M,n))\\
&\cong H^{n}\Mod_{A}(NL_{A/R}, M)\\
&\cong H^n\Der_R(X,M).
\end{align*}
Here $X \to A$ is a cofibrant replacement for $A$ as $R$-algebra.
\def\cH{{{\mathcal{H}}}}
In fact, we can define the Andr\'e-Quillen cohomology {\it space} by the equation
$$
\cH^n(A/R;M) = \map_{s\Alg_R/A}(X,A \ltimes K(M,n))
$$
where $X \to A$ is a cofibrant replacement. Then
$$
D^n(A/R;M) = \pi_0\cH^n(A/R;M) \cong \pi_t\cH^{n+t}(A/R;M).
$$
\end{remark}

\begin{remark}\label{applications}Andr\'e-Quillen (co-)homology arises naturally
in homotopy theoretic settings: see for example, \cite{MB},
\cite{bous-obstruct}, \cite{GJPAA}, and \cite{GH}. For example,
if $X$ and $Y$ are spaces and $\phi:X \to Y$ is a chosen basepoint
for the space of maps $\map(X,Y)$, then we leave it as an exercise to make
sense of the statement that there is a Hurewicz map
$$
\pi_t(\map(X,Y);\phi) \longr \Der_{\cK}(H^\ast Y, \Sigma^t H^\ast X).
$$
Here $H^\ast(-)$ is homology with $\FF_p$-coefficients, $\cK$ is the category
of unstable algebras over the Steenrod algebra, and $\Sigma^t H^\ast X$
is an $H^\ast Y$ module via $\phi^\ast$. This is the edge homomorphism  
of the Bousfield-Kan spectral sequence, and the $E_2$-term is given by
the higher derived functors of $\Der$ -- in short, by Andr\'e-Quillen cohomology.
\end{remark} 

\section{Resolutions in Model Categories}

In this section we turn to the project of creating and manipulating
simplicial resolutions in model categories. In previous sections we
worked with resolutions of modules and resolutions of algebras, but
now we'd like to turn to resolving objects which, at least for us, are
only inherently defined up to weak equivalence -- for example, we might
want to resolve topological spaces by spheres. This leads to the
resolution model categories of Dwyer, Kan, and Stover. We will
explain the approach developed by Bousfield.

\subsection{Skeletons and the skeletal decomposition}

\def\sk{{{\mathrm{sk}}}}

This preliminary section will be devoted to developing the canonical
and natural filtration of a simplicial object by its skeleton. This will
give us a chance to introduce some notation and prove, for example,
that every free morphism of simplicial algebras is a cofibration. 

The skeletons are defined by introducing some subcategories of
the ordinal number category $\ddelta$. Let $\ddelta[n] \subset \ddelta$
be the full-subcategory with objects $[k]$, $k \leq n$. If $\calC$ is
a category, we let $s_n\calC$ be the category of contravariant functors
from $\ddelta[n]$ to $\calC$. There is then a restriction functor
$r_\ast: s\calC \to s_n\calC$ from simplicial objects in $\calC$. If the
category $\calC$ has enough limits and colimits -- which we will assume --
this functor has a left adjoint, given by left Kan extension
$$
r^\ast: s_n\calC \longr s\calC.
$$
If $X \in s\calC$ is a simplicial object we define the $n$th {\it skeleton}
of $X$ by the formula
$$
\sk_nX = r^\ast r_\ast X.
$$
By the construction of Kan extensions, we have that
\begin{equation}\label{stage-m-skel}
(\sk_nX)_m = \underset{\phi:[m] \to [k]}\colim\ \phi^\ast X_k
\end{equation}
where $\phi$ runs over all the morphisms in $\ddelta$ with $k \leq n$.
But, since every such morphism factors uniquely as a surjection
followed by an injection, we may take the colimit over the {\it surjections}
in $\ddelta$ with $k \leq n$. In particular, $(\sk_nX)_m = X_m$
if $m \leq n$ and you should think of $\sk_nX$ as the sub-simplicial
set of $X$ generated by the non-degenerate simplices in degrees
$k \leq n$. There are natural maps $\sk_nX \to \sk_{n+1}X$
and $\sk_n X \to X$. Furthermore, the natural map
$$
\colim_n\ \sk_n X \longr X
$$
is an isomorphism.

More generally, if $X \to Y$ is a morphism of simplicial objects
in $\calC$, we define the relative $n$-skeleton $\sk_n^X Y$ by
the pushout diagram
$$
\xymatrix{
\sk_n X \rto \dto & X \dto\\
\sk_n Y \rto & \sk_n^X Y.
}
$$
Then $\sk_{-1}^X Y = X$ and $\colim_n \sk_n^X Y = Y$.

While these definitions are formal enough, they have impact because
we can explicitly describe the transition from the $(n-1)$st skeleton
to the $n$th skeleton. This is what we do next.

If $X$ is a simplicial set, we may define the degenerate $n$-simplices
as the image of the degeneracies:
$$
\bigcup_i s_i X_{n-1} \subseteq X_n.
$$
However, for more general
simplicial objects we need a categorical definition of the degeneracies.
As a hint, recall that the degeneracties in degree $n$ of a simplicial set
can equally be written
$$
\bigcup_{\phi:[n] \to [k]} \phi^\ast X_k = \colim_{\phi:[n] \to [k]}\ X_k
$$
where $\phi$ runs over the non-indentity surjections in $\ddelta$.

\begin{definition}\label{latching-matching} Let $\calC$ be a category
with all limits and colimits and let $X \in s\calC$ be a simplicial
object in $\calC$. Then the $n$th {\it latching object} of $X$ is
the colimit
$$
L_n X = \colim_{\phi:[n]\to [k]}\ X_k
$$
where $k$ runs over the non-identity surjections in the ordinal number
category $\ddelta$. The degeneracies define a natural map
$s:L_nX \to X_n$.

There is also a companion {\it matching object} given by the limit
$$
M_n X = \lim_{\psi:[k] \to [n]}\ X_k
$$
where $\psi$ runs over the non-identity injections. The face maps yield
a natural map $d:X_n \to M_nX$.
\end{definition}

By comparing Equation \ref{stage-m-skel} with the definition of the latching
object we immediately have a natural isomorphism
\begin{equation}\label{stage-n-skel-n-1}
(\sk_{n-1}X)_n = L_nX.
\end{equation}

Now we give $s\calC$ its standard simplicial structure, as in Equations
\ref{k-tensor-x} and following. A moment's thought should convince
the reader that if $Z \in \calC$ is regarded as a constant simplicial
object then there are natural isomorphisms
$$
s\calC(\Delta^n \otimes Z,Y) \cong \calC(Z,Y_n)
$$
and
$$
s\calC(\partial\Delta^n \otimes Z,Y) \cong \calC(Z,M_nY).
$$
If we set $Y = \sk_nX$ and $Z = X_n$ in the first equation, we obtain a natural map
$$
\Delta^n \otimes X_n \longr \sk_n X
$$
and if we set $Y = \sk_{n-1}X$ and $Z=X_n$ in the second equation,
we obtain a natural map
$$
\partial\Delta^n \otimes X_n \longr \sk_{n-1}X
$$
since $M_n\sk_{n-1}X = M_nX$. Setting $Z = L_nX$ and $Y = \sk_{n-1}X$
in the first equation we obtain a map
$$
\Delta^n \otimes L_nX \longr \sk_{n-1}X
$$
using Equation \ref{stage-n-skel-n-1}. In the end we obtain a commutative
diagram
$$
\xymatrix{
\Delta^n \otimes L_nX \sqcup_{\partial\Delta^n \otimes L_nX} \partial\Delta^n 
\otimes X_n \rto \dto & \sk_{n-1}X \dto\\
\Delta^n \otimes X_n \rto & \sk_nX.
}
$$
The claim, of course, is that this is a push-out. We can actually
give a relative version as well. See \cite{GJ} for an argument.

\begin{proposition}\label{skel-decomp}Let $X \to Y$ be a morphism
of simplicial objects and let
$$
L_n^XY = X_n \sqcup_{L_nX} L_nY.
$$
Then there is a push-out diagram
$$
\xymatrix{
\Delta^n \otimes L_n^XY \sqcup_{\partial\Delta^n \otimes L_n^XY} \partial\Delta^n 
\otimes Y_n \rto \dto & \sk_{n-1}^XY \dto\\
\Delta^n \otimes Y_n \rto & \sk_n^XY.
}
$$
\end{proposition}

Proposition \ref{skel-decomp} has a considerable simplification
if the morphism $X \to Y$ is $s$-free, as in Definition \ref{s-free}.

\begin{proposition}\label{s-free-skel} Suppose a morphism $X \to 
Y$ in $s\calC$ is $s$-free on a set of objects $Z_k$, $k \geq 0$.
Then there is a push-out diagram
$$
\xymatrix{
\partial\Delta^n \otimes Z_n \rto\dto & \sk_{n-1}^XY\dto\\
\Delta^n \otimes Z_n \rto & \sk_{n}^XY.
}
$$
\end{proposition}

This follows from Proposition \ref{skel-decomp} and the following observation: if $X$ is $s$-free on
$Z_k$, then
$$
L_n X \cong\underset{\phi:[n] \to [k]}\amalg\  \phi^\ast Z_k
$$
where $\phi$ runs over the non-identity surjections of $\ddelta$. Compare
Equation \ref{formula-for-k} -- this is the case of simplicial modules.

We can now prove the remaining implication of Proposition \ref{free-cof}.

\begin{lemma}\label{free-is-cob}Let $A \to B$ be a free map in
the category $s\Alg_R$. Then $A \to B$ is a cofibration.
\end{lemma}

\begin{proof} We need to show that any lifting problem
$$
\xymatrix{
A \rto \dto & X \dto^p\\
B \rto \ar@{-->}[ur] & Y
}
$$
with $p$ an acyclic fibration of simplicial $R$-algebras 
has a solution. We induct over the relative skeletons of $A \to B$
and solve the lifting problems
$$
\xymatrix{
\sk_{n-1}^AB \rto \dto & X \dto^p\\
\sk_n^AB \rto \ar@{-->}[ur] & Y
}
$$
By the definition of what it means to be free, the previous result
provides  a push-out square
$$
\xymatrix{
\partial\Delta^n \otimes S_R(P_n) \rto \dto & \sk_{n-1}^AB\dto\\
\Delta^n \otimes S_R(P_n) \rto & \sk_n^AB
}
$$
where $S_R(-)$ is the symmetric algebra functor and $P_n$ is a projective
$R$-module. Now $K \otimes S_R(P_n) = S_R(K \otimes P_n)$,
where $K \otimes (-)$ is interpreted in $R$-algebras or $R$-modules
as required. Thus, we need only solve the lifting problems
$$
\xymatrix{
S_R(\partial\Delta^n \otimes P_n) \rto \dto & X \dto^p\\
S_R(\Delta^n \otimes P_n) \rto \ar@{-->}[ur] & Y
}
$$
But this lifting problem in $s\Alg_R$ is adjoint to the lifting
problem
$$
\xymatrix{
\partial\Delta^n \otimes P_n \rto \dto_{i \otimes P_n} & X \dto^p\\
\Delta^n \otimes P_n \rto \ar@{-->}[ur] & Y
}
$$
in $s\Mod_R$, and this problem has a solution since $p$ remains
an acylic fibration in $s\Mod_R$ and $i \otimes P_n$ is a cofibration in
$s\Mod_R$.
\end{proof}

\subsection{Reedy model categories; the invariance of realization}

If $\calC$ has enough colimits, then there is a notion of geometric
realization from the category $s\calC$ of simplicial objects in $\calC$ back
down to $\calC$. In his thesis \cite{reedy}, Reedy noticed that if
$\calC$ is a model category, this could be made into a Quillen functor
and, in particular, realization preserved weak equivalences between
appropriately defined cofibrant objects. Since this is proved
using the technology of the previous section, we'll give an exposition
here. It also sets the stage for the resolution model categories of the
next section.

\begin{theorem}\label{reedy-mod}Let $\calC$ be a model category
and let $s\calC$ be the category of simplicial objects in $\calC$. Then
there is a model category structure on $s\calC$ where a morphism
$f:X \to Y$ is
\begin{enumerate}

\item a weak equivalence if $X_n \to Y_n$ is a weak equivalence
in $\calC$ for all $n \geq 0$;

\item a cofibration if the natural morphism from the relative
latching object
$$
L_n^XY = X_n \sqcup_{L_nX} L_n Y \to Y_n
$$
is a cofibration in $\calC$ for all $n\geq 0$; and

\item a fibration if the natural morphism to the relative matching
object
$$
X_n \to Y_n \times_{M_nY} M_nX
$$
is a fibration in $\calC$ for all $n \geq 0$.
\end{enumerate}
\end{theorem}

The proof is in \cite{reedy}, but see also \cite{GJ} and \cite{Local}.
We will refer to the weak equivalences, cofibrations, and fibrations
in this model category structure as {\it Reedy} weak equivalences, etc.

\begin{remark}\label{reedy-remark}1.) To get a feel for this
model category structure, note that if $X \in \calC$ is cofibrant,
then the constant object on $X \in s\calC$ is Reedy cofibrant.
However, constant objects are hardly ever Reedy fibrant.

2.) This model category structure
is not a simplicial model category structure in the standard 
simplicial structure on $s\calC$. (See Equation \ref{k-tensor-x} and
following.) If $i:K \to L$ is a cofibration of simplicial sets and
$j:X \to Y$ is a Reedy cofibration, then
$$
i \otimes j: K \otimes Y \sqcup_{K \otimes X} L  \otimes X \to
L \otimes Y
$$
is a Reedy cofibration which is a Reedy weak equivalence if $j$ is
a Reedy weak equivalence; however, it need not be a weak equivalence
if $i$ is a weak equivalence. Consider, for example, the inclusion
$$
d_0\otimes X :X =\Delta^0 \otimes X \to \Delta^1 \otimes X.
$$
However, if $\calC$ is a simplicial model category structure, then $s\calC$ inherits
an ``internal'' simplicial model category structure promoted up from
$\calC$.

3.) One example we might take for $\calC$ is the category $\sisets$
of simplicial sets. Then $s\calC$ is the category of bisimplicial
sets, but the Reedy model category structure priviledges the new
simplicial direction. Note that every bisimplicial set is automatically
Reedy cofibrant.

4.) At virtually the same time that Reedy was writing his thesis,
Bousfield and Kan wrote their book \cite{BK}. Central to the
existence of the homotopy spectral sequence of a cosimplicial
space was a model category structure on cosimplicial simplicial
sets. This is, it turns out, the Reedy structure,
after you've taken the opposite category enough times. One of their
theorems -- the homotopy invariance of the total space of a fibrant
cosimplicial space -- is essentially Theorem \ref{inv-of-realization}
below.

5.) If $\calC = \cgh$ is the category of compactly generated weak
Hausdorff spaces, the notion of a Reedy cofibrant object is a 
variation on the notion of a {\it proper} simplicial space. For proper,
one only requires that $L_n X \to X_n$ be a Hurewicz cofibration or,
perhaps, only a closed inclusion.
\end{remark}

Now suppose that $\calC$ itself is a simplicial model category. Let us
write 
$$
(-) \otimes_\calC (-):\sisets \times \calC \to \calC
$$
for the action of simplicial sets on $\calC$ -- reserving the symbol
$\otimes$ for the standard action of $\sisets$ on the category
$s\calC$ of simplicial objects. We will refer to these as the internal
and external actions respectively. 

We can use the internal action to define a geometric realization
functor for $s\calC$.

\begin{definition}\label{realize}Let $\calC$ be a simplicial model
category and let $X \in s\calC$ be a simplicial object in $\calC$.
Define the {\it realization} $|X| \in \calC$ of $X$ by the coequalizer
diagram in $\calC$.
$$
\xymatrix{
\coprod_{\phi:[m]\to [n]}\ \Delta^m \otimes_{\calC} X_n \ar@<.5ex>[r] \ar@<-.5ex>[r] & \coprod_n\  \Delta^n \otimes_\calC X_n
\rto & |X|.
}
$$
Here $\phi$ runs over the morphisms in $\ddelta$ and the parallel
arrows are obtained by evaluating on $\Delta^{(-)}$ and $X_{(-)}$
respectively.
\end{definition}

\def\diag{{{\mathrm{diag}}}}

\begin{examples}\label{real-exams}1.) This geometric realization functor from
simplicial sets to $\cgh$ of Definition \ref{geometric-realization}
can be obtained by regarding a simplicial set as a discrete simplicial
space and applying the geometric realization functor 
$|-|:s\cgh \to \cgh$.

2.) Here is an Eilenberg-Zilber Theorem.
If $X \in s\sisets$, then there is a natural isomorphism 
$|X| \to \diag(X)$ where $\diag(X)$ is the diagonal simplicial
set of the bisimplicial set $X$. To see this recall that the $k$-simplices
of $\Delta^n$ are the morphisms $\phi:[k] \to [n]$ in $\ddelta$.
If $X = \{ X_n\}$ is a simplicial object in simplicial sets, we can
define
$$
(\Delta^n\times X_n)_k = (\Delta^n)_k \times (X_n)_k \to
(X_k)_k
$$
sending $(\phi,y)$ to $\phi^\ast(y)$. It is an exercise to check
that there is an induced map of simplicial sets
$$
\coprod_n\  \Delta^n \times X_n \to \diag(X)
$$
and that this map coequalizes the diagram of Definition \ref{realize}.

3.) Let $A \in \calC$ and $K$ a simplicial set. Then using the external
action of $\sisets$ on $s\calC$ we get an object $K \otimes A \in
s\calC$. Then there is a natural isomorphism $|K \otimes A| \cong
K \otimes_{\calC} A$. The argument goes more or less as in the
previous example and is left as an exercise.
\end{examples}

The functor $|-|:s\calC \to \calC$ has a right adjoint
$$
A \mapsto \{ A^{\Delta^n}\}\ \defeq\ A^\Delta
$$
where $A^K$ is the exponential object internal to $\calC$. Then
we have the following variation on one of Reedy's main results.

\begin{theorem}\label{inv-of-realization} Let $\calC$ be a
simplicial model category. The adjoint pair
$$
\xymatrix{
|-|:s\calC \ar@<.5ex>[r] & \ar@<.5ex>[l] \calC: (-)^\Delta
}
$$
is a Quillen functor from $s\calC$ to $\calC$. In particular,
realization preserves cofibrations and weak equivalences between
cofibrant objects.
\end{theorem}

Rather than prove this result, let us examine how the skeleton
filtration interacts with realization. If $X \in s\calC$ 
there is an isomorphism $\colim |\sk_nX| \cong |X|$ and, using
Proposition \ref{skel-decomp} and Example \ref{real-exams}.3 we have
a push-out diagram
$$
\xymatrix{
\Delta^n \otimes_\calC L_nX \sqcup_{\partial\Delta^n \otimes_\calC L_nX} \partial\Delta^n 
\otimes_\calC X_n \rto \dto &| \sk_{n-1}X| \dto\\
\Delta^n \otimes_\calC Y_n \rto & |\sk_nX|.
}
$$
If $X$ is Reedy cofibrant, then the morphisms $L_nX \to X_n$ are
cofibrations, and then the corner axiom for simplicial model categories
(See Proposition \ref{other-sm7s}) implies that $|\sk_{n-1}X|
\to |\sk_{n}X|$ is a cofibration. 

In the case when $X$ is $s$-free -- which happens surprisingly
often in applications -- this push-out diagram specializes
even more to a diagram that really looks like something
you'd have with a simplicial complex. If $X$ is $s$-free
on $\{Z_k\}$, then Proposition \ref{s-free-skel} yields
a push-out diagram
$$
\xymatrix{
\partial\Delta^n \otimes_\calC Z_n \rto\dto & |sk_{n-1}X| \dto\\
\Delta^n \otimes_\calC Z_n \rto & |sk_{n}X|.
}
$$

To specialize to an application, if $X \in s\cgh$ is a simplicial space,
then the cofiber of $|\sk_{n-1}X| \to |\sk_{n}X|$ is the pointed 
space
$$
|\Delta^n/\partial\Delta^n| \wedge (X_n/L_nX) \cong 
S^n \wedge (X_n/L_nX)
$$
and we check that if $E_\ast$ is any generalized homology theory,
then
$$
\tilde{E}_{p+q}(|\Delta^p/\partial\Delta^p| \wedge (X_p/L_pX))
\cong N\tilde{E}_qX_p.
$$
Thus it follows that there is a spectral sequence, for any Reedy cofibrant
simplicial space $X$
\begin{equation}\label{homology-ss}
\pi_pE_q(X) \Longrightarrow E_{p+q}|X|.
\end{equation}
A variation of these ideas, using cosimplicial spaces
and homotopy yields the extremely useful Bousfield-Kan
spectral sequence of \cite{BK}, which is the prototype of all
second-quadrant homotopy spectral sequences.

\subsection{Resolution model categories}

\def\Ho{{{\mathbf{Ho}}}}
\def\cS{{\calC}}
\def\cP{{{\mathcal{P}}}}

The Reedy model category structure of the previous section
has two drawbacks. First, our simplicial
objects are often built as resolutions and, in so doing, we may want
to specify the projective objects from which we build the resolutions.
The Reedy model category structure doesn't do this -- there are
too many cofibrant objects. Second, after taking a resolution, we then
would want to examine its realization, typically through some spectral
sequence such as the homology spectral sequence of 
Equation \ref{homology-ss}. While a Reedy weak equivalence will
certainly yield an isomorphism of spectral sequences, there are many
other morphisms which also will yield an isomorphism of spectral sequences --
there are not enough weak equivalences in $s\calC$. We remedy both
problems as once. The basic ideas go back to Dwyer, Kan, and Stover
in \cite{DKS1} and \cite{DKS2}, but they were greatly expanded in
\cite{jardine-res} and, especially, \cite{bous-stover}.\footnote{Bousfield's
paper is written cosimplicially, rather than simplicially, but the arguments
are so categorical that they readily translate to the simplicial setting.}

For convenience and simplicity of language, we will work in a {\it pointed}
model category $\calC$ -- that is, a category where the unique map
from the initial object to terminal object is an isomorphism. Thus,
for example, we will work with pointed spaces $\cgh_\ast$.
We will also assume we
have a simplicial model category. Both assumptions are unnecessary,
as explained in \cite{bous-stover}, and can be removed at the cost of
ramping up the language.

If $A \in \calC$ is cofibrant, the simplicial structure gives a canonical
model for $\Sigma A$, the suspension of $A$ as the cofiber of
$$
\partial\Delta^1 \otimes A \sqcup_{\partial\Delta^1 \otimes\ast}
\Delta^1 \otimes \ast \longr \Delta^1 \otimes A.
$$
Here $\ast$ is the initial (equals terminal) object. The cone,
$\Cone(A)$, on $A$ has a similar description.

We begin by specifying the building blocks of our resolutions.

\begin{definition}\label{homotopy-cogroup} Let $\calC$ be our fixed
simplicial model category and let $\Ho(\calC)$ denote the homotopy
category. Then a {\it homotopy cogroup object} $A \in \calC$ is a
cofibrant object
so that $\Ho(\calC)(A,-)$ is a functor to groups. A set of {\it projectives}
in $\calC$ will be a set of homotopy cogroup objects which is
closed under finite coproducts and suspension.
\end{definition}

\begin{examples}\label{exams-cogroup}1. In $\cgh_\ast$, any
suspension is a homotopy cogroup object. Indeed, if $X$ is
any pointed CW-complex, we get a set of projectives
by taking finite wedges of the spaces $\Sigma^n X$, $n \geq 1$.
In particular, the example of \cite{DKS1} is given by setting
$X = S^0$, whence the set of projectives is finite wedges of
positive spheres.

2. A rich class of examples arises in spectra
by considering homology theories for which there is
a universal coefficient spectral sequence. See \cite{GH}. Spectra
have not been discussed here, but there is a very rich literature.
For a relatively straightforward model category on spectra whose
homotopy category is the stable category, see \cite{BF}.

3. Let $\calC = s\Mod_R$ for a commutative ring $R$. Then any
cofibrant object is a homotopy cogroup object. An obvious
set of projectives is given by taking finite sums of objects of the
form $\Sigma^n R$ where $R$ is regarded as a constant simplicial
$R$-module. Note $\Sigma^n R \cong K(R,n)$.

4. The category $s\Alg_R$ is not pointed -- and this is an example
where we might want the unpointed generalization of this discussion
here. In any case, if $M \in s\Mod_R$ is any simplicial cofibrant
simplicial $R$-module,
then the symmetric algebra $S_R(M) \in s\Alg_R$ is a homotopy
cogroup object, in the obvious sense, and we might want to consider
the set of projectives generated by $S_R(\Sigma^n R)$.
\end{examples}

The following ideas are meant to echo the idea of projectives
and surjections introduced on the very first page of these notes.

\begin{definition}\label{projectives} Let $\calC$ be a pointed simplicial model category and $\Ho(\calC)$ its homotopy category. Fix a set $\cP$ of
projectives in $\calC$. Write $[-,-]$ for $\Ho(\calC)(-,-)$.
\begin{enumerate}

\item A morphism $p: X \to Y$ in $\Ho(\calC)$ is
$\cP$-{\it epi} if $p_\ast:[P,X] \to [P,Y]$ is onto for each $P \in \cP$.

\item An object $A \in \Ho(\calC)$ is $\cP$-{\it projective} if
$$
p_\ast:[A,X] \longr [A,Y]
$$
is onto for all $\cP$-epi maps.

\item A morphism $A \to B$ in $\calC$ is called a $\cP$-{\it projective
cofibration} if it has the left lifting property for all $\cP$-epi
fibrations in $\calC$.
\end{enumerate}
\end{definition}

The classes of $\cP$-epi maps and of $\cP$-projective objects determine
each other; furthermore, every object in $\cP$ is $\cP$-projective.
Note however, that the class of $\cP$-projectives is closed under
arbitrary coproducts. The class of $\cP$-projective cofibrations will
be characterized below; see Lemma \ref{proj-cofib}.

\begin{lemma}\label{enough-projective}The category $\Ho(\cS)$ has
enough $\cP$-projectives; that is, for every object $X \in \Ho(\cS)$
there is a $\cP$-epi $Y \to X$ with $Y$ $\cP$-projective.
\end{lemma}

\begin{proof} We can take 
$$
Y = {\underset{P \in \cP}{\amalg}}\ {\underset{f:P\to X}{\amalg}} P
$$
where $f$ ranges over all maps $P \to X$ in $\Ho(\cS)$.
\end{proof}

We now come to the $\cP$-resolution model category structure. Recall
that a morphism $f:A \to B$ of simplicial abelian groups is a
weak equivalence if $f_\ast:\pi_\ast A \to \pi_\ast B$ is an
isomorphism. Also $f:A \to B$ is a fibration if the induced map
of normalized chain complexes $Nf:NA \to NB$ is surjective in
positive degrees. We haven't been explicit about this, but the
same statements hold for simplicial groups. If $A$ is a simplicial
group, not necessarily abelian, we can still form the normalized
complex $NA$ as in Example \ref{moore-complex}.
We easily check that a morphism
$f:A \to B$ of simplicial groups is a fibration if and only
if it $Nf$ is surjective in positive degrees.

\begin{definition}\label{P-structures} Let $\calC$ be a pointed simplicial
model category with a fixed set of projectives $\cP$. Let
$f:X \to Y$ be a morphism in the category $s\calC$ of simplicial
objects in $\calC$. Then
\begin{enumerate}

\item the map $f$ is a $\cP$-{\it equivalence} if the induced morphism
$$
f_\ast:[P,X] \longr [P,Y]
$$
is a weak equivalence of simplicial groups for all $P \in \cP$; 

\item the map $f$ is a $\cP$-{\it fibration} if it is a Reedy
fibration and $f_\ast:[P,X] \longr [P,Y]$ is a fibration of simplicial
groups for all $P \in \cP$;

\item the map $f$ is a $\cP$-{\it cofibration} if the induced maps
$$
X_n \sqcup_{L_nX} L_nY \longr Y_n, \qquad n \geq 0,
$$
are $\cP$-projective cofibrations.
\end{enumerate}
\end{definition}

Then, of course, the theorem is as follows. 

\begin{theorem}\label{res-model} With these definitions of $\cP$-equivalence,
$\cP$-fibration, and $\cP$-cofibration, the category $s\calC$ becomes a
simplicial model category.
\end{theorem}

The proof is given in \cite{bous-stover}. We call this the
$\cP$-resolution model category structure.  It is cofibrantly
generated if $\calC$ is cofibrantly generated and an object is
$\cP$-fibrant if and only if it is Reedy
fibrant. Furthermore, any Reedy weak equivalence
is a $\cP$-equivalence and any $\cP$-cofibration is a Reedy
cfoibration. The next result gives a characterization of $\cP$-cofibrations.

Define a morphism $X \to Y$ in the category $\calC$ to be $\cP$-{\it free}
if it can be written as a composition 
$$
\xymatrix{
X \rto^i & X \amalg F  \rto^q &Y\\
}
$$
where $i$ is the inclusion of the summand, $F$ is cofibrant and
$\cP$-projective, and $q$ is an acyclic cofibration. The following
is also in \cite{bous-stover}.

\begin{lemma}\label{proj-cofib} A morphism $X \to Y$ in $\calC$ is a
$\cP$-projective cofibration if and only if it is a retract of $\cP$-free
map.
\end{lemma}

\def\Cone{{{\mathrm{Cone}}}}

\begin{remark}\label{stover-resolution}In his paper that inaugurated
this subject \cite{stover}, Stover wrote
down a very concrete model for the $\cP$-cofibrant replacement in
$s\calC$, at least in the case where every object is fibrant. If $P \in \cP$,
let $\Cone(P)$ denote the cone on $P$. For
$A$ in $\calC$, define $V(A)$ by the push-out diagram
$$
\xymatrix{
\amalg_P \amalg_{\Cone(P) \to A} P \rto \dto & 
\amalg_P \amalg_{P \to A} P \dto \\
\amalg_P \amalg_{\Cone(P) \to A} \Cone(P) \rto & V(A).\\
}
$$
where the coproducts are over morphisms in $\calC$.
There are natural maps $V(A) \to A$ and $V(A) \to V^2(A)$. The first
is given by evaluation and the second by noticing that each summand
in the pushout diagram defines a morphism $P \to V(A)$ (or
$\Cone(P) \to V(A)$). Thus,
by iterating $V$, we obtain an augmented simplicial object 
$V_\bullet(A) \to A$. Stover proves that this is a $\cP$-equivalence
and that $V_\bullet(A)$ is $\cP$-cofibrant. It is also easy to prove that there is a (non-natural) 
homotopy equivalence $V(A) \to W$, where $W$ is isomorphic
to a coproduct of projectives in $\cP$.
\end{remark}

\begin{examples}1. In the case where $\calC = \cgh_\ast$ and
$\cP$ is the collection of finite wedges of spheres, a cofibrant
replacement $X \to A$ for a space $X$ is a ``resolution'' by
spheres. The map $|X| \to A$ will be a weak equivalence if
$A$ is path-connected or, more generally, will be a weak
equivalence from $|X|$ to the component of the basepoint of
$A$. These simplicial resolutions were used by Stover to
define a spectral sequence for studying the homotopy groups
of a wedge.

2. If we take $\cP$ to be the collection of finite wedges of
spheres $S^k$, $k \geq n$, then $|X| \to A$ will be a model
for the $(n-1)st$ connected cover of $A$.

3. In \cite{GH}, we resolved spectra by finite CW spectra $P$ with the
property that $E_\ast P$ was projective for some homology theory
$E$. This was an enormous help in developing a spectral sequence for
studying the homotopy type of mapping spaces of structured ring
spectra.

4. One big impact of these resolutions is to realization problems.
For example, if $A$ is a pointed space, then $\pi_\ast A$ has a great
deal of structure coming from Whitehead products, compositions, etc.
Writing down axioms for this structure yields the notion of a 
$\Pi$-algebra and one can ask: given a $\Pi$-algebra $\Lambda$,
in there a space $A$ so that $\pi_\ast A \cong \Lambda$? In
\cite{BDG} we defined and studied the space of {\it all} such realizations,
and gave a decomposition of that space as an inverse limit of
a tower where the layers are governed by Andr\'e-Quillen cohomology.
A transportation of these ideas to \cite{GH} yielded solutions to
similar moduli problems in stable homotopy theory. Indeed, one
statement of the celebrated Hopkins-Miller theorem is that the
moduli space of all $E_\infty$-ring structures on Morava's $E_n$
spectrum is a space of the form $BG$, where $G$ is a group -- in fact,
the Morava stabilizer group. The basic algebraic fact at work
is that $(E_n)_\ast \to (E_n)_\ast E_n$ is formally \'etale and,
hence, has vanishing Andr\'e-Quillen cohomology.
\end{examples}

\begin{example}\label{ss-simp-model} We meet a precursor of 
resolution model categories in Quillen's original work.
See \cite{HA} \S II.6. If
$M$ and $N$ are two simplicial $R$-modules, let us define
$M \otimes^L_R N$ to be the total left derived functor of
$M \otimes_R (-)$ applied to $N$; that is, up to weak equivalence
$M \otimes_R^L N = M \otimes_R^L Q$ where $Q$
is a cofibrant replacement for
$N$. In \cite{HA}, Quillen writes down a spectral sequence
$$
\mathrm{Tor}^R_p(\pi_\ast M, \pi_\ast N)_q \Longrightarrow \pi_{p+q}(M \otimes_R^L N).
$$
The grading $q$ arises from the grading on the homotopy groups.
Quillen's argument is essentially as follows.

Consider the resolution model category structure on $s(s\Mod_R)$ defined
by the homotopy cogroup objects $\Sigma^n R = K(R,n)$. See
Example \ref{exams-cogroup}.3. Then if we
regard $N \in s\Mod_R$ as a constant simplicial object in $s(s\Mod_R)$,
we can form a $\mathcal{P}$-cofibrant replacement $P_\bullet \to N$
for $N$. Taking homotopy in the old (internal) simplicial degree first,
then in the new (external) simplicial degree, we have that
$\pi_0\pi_\ast P_\bullet \cong \pi_\ast N$
and $\pi_p\pi_\ast P_\bullet = 0$. Furthermore, $P_p$ is a cofibrant
simplicial module homotopy
equivalent to a direct sum of simplicial modules of the form
$K(R,n)$. Together, these two facts imply that $Q = \mathrm{diag}P_\bullet
\to N$ is a cofibrant replacement.

To get the spectral sequence, form the bisimplicial module $M \otimes_R
P_\bullet = \{M \otimes_R P_p\}$ and filter by the new (external)
simplicial degree $p$. Then the
$E^1$ term is given by
$$
E^1_{p,q} = \pi_p (M \otimes_R \pi_\ast P_p)
$$
and we have, from the previous paragraph, that $\pi_\ast P_\bullet \to \pi_\ast
N$ is a graded projective resolution. So Quillen's spectral sequence
follows.
\end{example}

\input biblio
\end{document}

%% file: macros.tex
\def\cgh{{{\mathbf{CGH}}}}
\def\Ch{{{\mathbf{Ch}}}}
\def\calC{{{\mathcal{C}}}}
\def\calD{{{\mathcal{D}}}}
\def\cO{{{\mathcal{O}}}}
\def\longr{{\longrightarrow}}
\def\ZZ{{{\mathbb{Z}}}}
\def\FF{{{\mathbb{F}}}}
\def\cK{{{\mathcal{K}}}}
\def\Der{{{\mathrm{Der}}}}
\def\Cone{{{\mathrm{Cone}}}}
\def\sisets{{{s\mathbf{Sets}}}}

\def\Ab{{{\mathbf{Ab}}}}
\def\Der{{{\mathrm{Der}}}}
\def\colim{{{\mathrm{colim}}}}

\def\Mod{{{\mathbf{Mod}}}}
\def\defeq{{{\overset{\mathrm{def}}=}}}
\def\ddelta{{{{\Delta}}}}